%% file: SIAM_article.tex
\documentclass[onefignum,onetabnum]{siamart171218}

\input{ex_shared}

\ifpdf
\hypersetup{
  pdftitle={Unified Learning of the Profile Function in Discrete Keller-Segel Models},
  pdfauthor={Chi-An Chen, Chun Liu, Ming Zhong}
}
\fi

\usepackage{bm}
\DeclareMathOperator*{\argmin}{arg\,min}
\newcommand{\f}{\mathbf{f}}
\newcommand{\F}{\mathbf{F}}
\newcommand{\br}{\mathbf{r}}
\newcommand{\bu}{\mathbf{u}}

\newcommand{\x}{\mathbf{x}}
\newcommand{\X}{\mathbf{X}}
\newcommand{\y}{\mathbf{y}}
\newcommand{\Y}{\mathbf{Y}}
\DeclareMathOperator{\dif}{d \!}        
\newcommand{\dt}{\dif t}
\newcommand{\dx}{\dif\x}
\newcommand{\zero}{\mathbf{0}}
\newcommand{\E}{\mathbb{E}}
\newcommand{\R}{\mathbb{R}}
\newcommand{\norm}[1]{| #1 |}

\usepackage{dcolumn}

\begin{document}

\maketitle

\title{Learning of the Profile Function in Deterministic and Stochastic Keller-Segel Models}

\date{\today}

\begin{abstract}
We propose a unified learning framework for identifying the profile function in discrete Keller-Segel equations, which are widely used mathematical models for understanding chemotaxis. Training data are obtained via either a rigorously developed particle method designed for stable simulation of high-dimensional Keller–Segel systems, or stochastic differential equations approximating the continuous Keller-Segel PDE. Our approach addresses key challenges, including data instability in dimensions higher than two and the accurate capture of singular behavior in the profile function. Additionally, we introduce an adaptive learning strategy to enhance performance. Extensive numerical experiments are presented to validate the effectiveness of our method.
\end{abstract}
\maketitle
%
\begin{center}
First circulated and This version: October $21$, $2025$ 
\end{center}
\section{Introduction \label{sec:introduction}}
Chemotaxis, the directed movement of cells or organisms in response to chemical gradients, plays a fundamental role in a wide range of biological processes, including embryonic development, immune system responses, and tumor progression~\cite{painter2018mathematicalmodelschemotaxisapplications, chavanis2010stochastic, zheng2023stability, rocca2023cahn}. The so-called Keller-Segel model provides a fundamental mathematical framework to describe this behavior, which couples the dynamics of cell density with that of a chemoattractant concentration~\cite{keller1971model}. These equations are known for their rich and sometimes singular behavior, including the emergence of self-organization, pattern formation, and finite-time blow-up, critically depending on the nature of the interaction profile~\cite{jin2016pattern, bedrossian2011local, raphael2014stability, bian2013dynamic, cao2025finite, collot2023collapsing}.

The Keller-Segel model has been extensively studied both analytically and numerically in order to better understand the mechanisms underlying chemotactic aggregation~\cite{blanchet2006two, arumugam2021keller, huang2017error, bresch2019mean, wang2024deepparticle, benzakour2020linearized, carrillo2019blob, barwolff2018numerical, bedrossian2011local, jin2016pattern, raphael2014stability, bian2013dynamic, patlak1953random}. A central challenge in these efforts lies in accurately characterizing the interaction profile, particularly in high-dimensional or singular regimes where traditional modeling assumptions may break down. Conventional approaches often rely on analytically prescribed forms for the interaction profile, which may fail to capture the complexity of observed data or microscale behaviors. Recent advances in machine learning have introduced new paradigms for modeling such systems, particularly through data-driven frameworks that allow the inference of governing equations or functional components directly from observed trajectories~\cite{wang2024deepparticle, benzakour2020linearized, hwang2021neural}. In this work, we investigate a unified variational learning approach originally developed for interacting particle systems~\cite{LZTM2019, ZMM2020, MMQZ2021, FMMZ2022, MTZM2022, feng2023learning}, and adapt it to recover the interaction profile from simulation data generated by the Keller–Segel equations. Our objective is to learn a non-parametric, data-driven representation of the interaction function that best captures the observed dynamics, thereby enhancing both interpretability and predictive performance in chemotaxis modeling.

A critical component of our methodology is a particle-based numerical solver tailored for high-dimensional Keller–Segel systems~\cite{wang2024deepparticle, carrillo2019blob, craig2016blob, chertock2006practical, chertock2001particle}. This method enables efficient and robust simulation in dimensions higher than three, where traditional grid-based PDE solvers often become unstable or computationally prohibitive. Furthermore, the particle method provides a natural means for resolving fine-scale or singular behavior—such as sharp density gradients or finite-time blow-up—that are essential features in chemotactic systems~\cite{chertock2024hybrid}. Using this framework, we generate high-fidelity datasets that serve as the foundation for our variational learning pipeline. To further improve learning accuracy and generalization, we incorporate an adaptive refinement strategy that dynamically adjusts the resolution of the function approximation in response to observed discrepancies. This adaptive scheme enables localization of regions with high curvature or emergent singularities, where finer approximations are crucial. Furthermore, we extend the particle method to simulate a two-dimensional stochastic Keller–Segel model with weak diffusion, highlighting its applicability to simulate stochastic particle systems. This also connects to a growing body of research on learning particle dynamics from stochastic differential equations~\cite{lu2020learning, guo2024learning, guo2024noiseguidedstructurallearning}. We demonstrate the effectiveness of our approach through a series of numerical experiments conducted across multiple spatial dimensions and parameter regimes, illustrating the method’s ability to recover accurate and physically meaningful interaction profiles.

The remainder of this paper is organized as follows. In Section~\ref{sec:model}, we introduce the energy-based formulation of the Keller–Segel model and describe the deterministic particle method as well as the stochastic particle trajectories governed by SDEs used to simulate the dynamics. Section~\ref{sec:learn} presents our updated variational learning framework and the adaptive refinement strategy, applicable to both deterministic and stochastic particle data. Numerical results and validation studies are provided in Section~\ref{sec:examples}, demonstrating the efficacy of our approach. Finally, we conclude in Section~\ref{sec:conclude} with a discussion of potential extensions and future directions.
\section{Model Description\label{sec:model}}
In this section, we discuss the general Keller–Segel model:
\begin{equation}\label{Keller_Segel_equation}
\partial_t \rho = \Delta \rho + \chi \nabla \cdot \bigl(\rho \, \nabla (W * \rho)\bigr),
\end{equation}
where $\rho = \rho(t, x)$ with $(t, x) \in [0, T]\times\Omega, \Omega \subset \mathbb{R}^d$. Here $\rho$ denotes the density of the species under consideration, $\chi$ is the sensitivity parameter describing the response of the bacteria to the chemoattractant, and W is a kernel function governing the attractive interaction. We derive this PDE formulation from an energy–dissipation law for a general aggregation–diffusion energy model, and we also discuss the stochastic-diffusion perspective of the Keller–Segel model. These derivations will aid in the design of our learning framework later on.
\subsection{Aggregation Diffusion Model\label{sec:aggregation_diffusion}}
An aggregation–diffusion model captures the evolution of a density under the competing effects of dispersal and attraction. The diffusion term represents the tendency of particles or individuals to spread out due to random motion, while the aggregation term accounts for nonlocal interactions that drive particles toward each other. This class of models appears in many contexts — such as chemotaxis, biological swarming, granular media, and gravitational systems — and captures a wide range of phenomena, including spreading, pattern formation, and finite-time concentration. Such models typically satisfy a general energy–dissipation law: their dynamics are governed by an energy functional and a dissipation mechanism, while the continuity equation ensures mass conservation:
\begin{equation} \label{energy-dissipation}
\left\{
\begin{aligned}
& \frac{\dif}{\dt} \mathcal{E}[\rho](t) = -\mathcal{D}[\rho](t), \\
& \partial_t \rho + \nabla \cdot (\rho\bu) = 0, \quad \bu = -\nabla \left( \frac{\delta \mathcal{E}}{\delta \rho} \right).
\end{aligned}
\right.
\end{equation}
Here, $\rho$ denotes the density of the system, $\mathcal{E}$ the free energy, $\mathcal{D}$ the dissipation functional, $\mathbf{u}$ the velocity field of the system, and $\frac{\delta \mathcal{E}}{\delta \rho}$ the variational derivative of the free energy with respect to the density $\rho$. The dissipation functional $\mathcal{D}$ is given by
\begin{equation}\label{dissipation}
\mathcal{D}[\rho] = \int_{\R^d} \rho |\bu|^{2} \, \dx.
\end{equation} 
The energy functional $\mathcal{E}$ typically takes the form
\begin{equation}\label{energy}
\mathcal{E}[\rho](t) = \int_{\R^d} \big(\rho(t, \x)V(\x) + \frac{1}{2} (W(\x) * \rho(t, \x)) \rho(t, \x) + \rho F_{m}(\rho)(t, \x)\big) \, \dx,
\end{equation}
where
\begin{equation}
F_{m}(\rho) =
\begin{cases}
\log \rho, & m = 1, \\
\frac{\rho^{m-1}}{m-1}, & m > 1.
\end{cases}
\end{equation}
Here, $V$ is the external potential representing environmental effects, $W$ is a radially symmetric interaction kernel; and the last term corresponds to the diffusion contribution, and $*$ denotes the usual convolution operator. The parameter $m$ characterizes different diffusion regimes: for $m = 1$, the model reduces to the classical heat equation, describing linear diffusion where the spreading rate is independent of the solution’s magnitude; for $m > 1$, it models nonlinear diffusion, as in the porous medium equation, where diffusion slows in low-density regions, relevant to phenomena such as gas flow in porous media or population dynamics.  The energetic variational approach combines the Least Action Principle applied to the free energy with the Maximum Dissipation Principle applied to the dissipation energy \cite{giga2017variational}. To derive the partial differential equations from \eqref{energy-dissipation}, we introduce the flow map $\mathbf{x}_{E}(t, \mathbf{x}_{L})$ with $\mathbf{x}_{E}(t,\cdot):\mathbb{R}^{d}\to\mathbb{R}^{d}$, where $\mathbf{x}_{E}$ denotes the Eulerian coordinate and $\mathbf{x}_{L}$ the Lagrangian coordinate. This flow map is associated with the velocity field $\mathbf{u}(\mathbf{x}_{E},t)$. For a given $\mathbf{x}_{L}$, the flow map satisfies
\begin{equation}
\dot{\x}_{E} = \frac{\dif}{\dt} \x_{E}(t, \x_{L}) = \mathbf{u}(t, \x_{E}(t, \x_{L})), \quad \x_{E}(0, \x_{L}) = \x_{L}.
\end{equation}
We will be focusing on the Eulerian coordinates from now, thus by default we use $\mathbf{x}_E = \mathbf{x}$, unless defined otherwise.  By a standard energetic variation procedure (see \cite{giga2017variational} for details), we obtain the force balance equation
\begin{equation}    
\frac{\delta \mathcal{E}}{\delta \x} + \frac{\delta \mathcal{D}}{\delta \dot{\mathbf{x}}} = 0.
\end{equation}
For the specific choices of the energy $\mathcal{E}$ in \eqref{energy} and the dissipation $\mathcal{D}$ in \eqref{dissipation}, together with the continuity equation in \eqref{energy-dissipation}, this variational principle leads to the PDE
\begin{equation}
\partial_{t}\rho(t, \x) = \nabla \cdot (\rho(t, \x) \nabla V(\x) + \nabla \cdot ((\nabla W(\x) * \rho(t, \x)) \rho(t, \x) + \Delta (\rho(t, \x))^{m}.
\end{equation}
\subsection{A Deterministic Particle Method for the Keller–Segel Model\label{sec:deterministic}}
To develop a numerical approximation of the force balance equation, we adopt a discrete variational approach by first discretizing the energy-dissipation law \eqref{energy-dissipation} at the particle level. Specifically, we approximate the initial measure $\rho_{0}$ as a finite sum of $N$ Dirac delta functions centered at particle positions, leading to a particle representation of the free energy and dissipation. By taking variations of the discrete free energy and dissipation with respect to particle trajectories, we obtain a system of ODEs governing the particle dynamics. In the absence of the diffusion term in \eqref{energy-dissipation}, the deterministic particle method follows this discrete variational framework and represents the density as
\[
\rho_{0}(\x) \approx \rho_{0}^{N}(\x) = \sum_{i=1}^{N} \omega_{i}\delta(\x - \x_{i}), \quad \x_{i} \in \mathbb{R}^{d}, \quad \omega_{i} \geq 0,
\]
where $\omega_{i}$ denotes the mass assigned to the i-th particle.
Under suitable regularity on the potential $V$ and the interaction kernel $W$, the solution $\rho^{N}$ with initial data $\rho_{0}^{N}$ retains its form as a sum of Dirac delta functions for all time $t$. This can be expressed as:  
\begin{equation}\label{discrete_density}
    \rho^{N}(t,\x) = \sum_{i=1}^{N}\omega_{i}\delta(\x - \x_{i}(t)).
\end{equation}
However, when the diffusion term is present, we cannot evaluate $\mathcal{F}^{m}(\rho)$ directly on a sum of Dirac delta functions. This issue arises because the diffusion term $\mathcal{F}^{m}(\rho)$ involves nonlinear functionals of $\rho$, which are not well-defined for singular measures. To address this, we apply a regularization only to the density appearing in the diffusion term by convolving it with a smoothing kernel \cite{carrillo2019blob, craig2016blob}. This leads to the following regularized expression:
\begin{equation}
    \mathcal{F}^{m}_{h}(\rho) = \int \rho F_{m}(K_{h} * \rho) \ \dx.
\end{equation}
Here, $K_{h}(\cdot)$ is a mollifier, commonly chosen as a Gaussian kernel. Importantly, this regularization is applied only within the diffusion term, while the other terms in the energy-dissipation law remain unchanged. As a result, the particle representation of the density is preserved throughout the evolution, ensuring that the discretization maintains the underlying particle structure. Moreover, we use the regularized free energy $\mathcal{E}_{h}(\rho)$ as an approximation to the original free energy $\mathcal{E}$, incorporating the regularized diffusion term $\mathcal{F}^{m}_{h}(\rho)$ replacing $\F^m$ in \eqref{energy}. When we substitute \eqref{discrete_density} into the regularized free energy, we obtain 
\begin{equation*}
\mathcal{E}_{h}(\{\x_{i}(t)\}_{i=1}^{N}) = \sum_{i=1}^{N}\omega_{i}\big(V(\x_{i}(t)) + \frac{1}{2}\sum_{j=1}^{N}\omega_{j} W(\mathbf{r}_{i, j}(t)) + F_{m}(\sum_{j=1}^{N}\omega_{j}K_{h}(\mathbf{r}_{i,j}(t)))\big),
\end{equation*}
where $\mathbf{r}_{i, j}(t) = \x_i(t) - \x_j(t)$.  Similarly, substituting \eqref{discrete_density} into the dissipation term yields 
\begin{equation}
    \mathcal{D}_{h}(\{\x_{i}\}_{i=1}^{N}; \{\dot{\x}_{i}\}_{i=1}^{N}) = \frac{1}{2} \sum_{i=1}^{N} \omega_{i} |\dot{\x}_{i}|^{2}.  
\end{equation}
By taking the variation with respect to each particle $\x_{k}(t), \ k = 1,\cdots, N$, we obtain 
\[
\begin{aligned}
    \frac{\delta \mathcal{E}_{h}}{\delta \x_{k}}  & = \omega_{k} \Big[\nabla V(\x_{k}) + \sum_{\substack{i=1 \\ i\neq k}}^{N} \omega_{i} \nabla W(\br_{k, i}) 
      + \sum_{\substack{i=1 \\ i\neq k}}^{N} \omega_{i} F_{m}'\left(\sum_{j=1}^{N}\omega_{j} K_{h}(\br_{k, j})\right) \nabla K_{h}(\br_{k, i}) \Big]\\
     & + F_{m}'\sum_{j=1}^{N}\omega_{j}K_{h}(\br_{k, j})\sum_{i=1}^{N} \omega_{i}\nabla K_{h}(\br_{k, i}), \quad \br_{k, i} = \x_k - \x_i.
\end{aligned}
\]
For the dissipation term, we take the variation with respect to $\dot{\x}_{k}$ for $k = 1,\cdots, N$: 
\begin{equation}
    \frac{\delta \mathcal{D}_{h}}{\delta \dot{\x}_{k}} = \omega_{k} \dot{\x}_{k}.
\end{equation}
According to the force balance, we have 
\begin{equation}
    \frac{\delta \mathcal{D}_{h}}{\delta \dot{\x}_{k}} + \frac{\delta \mathcal{E}_{h}}{\delta \x_{k}} = 0,
\end{equation}
which gives the ODE system of the particle trajectories  
\begin{equation}\label{ODE_sys}
    \dot{\x}_{k}(t) = -\frac{1}{\omega_{k}} \frac{\delta \mathcal{E}_{h}}{\delta \x_{k}}(\{\x_{i}(t)\}_{i=1}^{N})
\end{equation}
Using the fully implicit Euler scheme, we can discretize \eqref{ODE_sys} as follows:
\begin{equation}\label{particle_evolution_deter}
\begin{split}
    \frac{\x^{n+1}_{k} - \x^{n}_{k}}{\tau} = -\frac{1}{\omega_{k}} \frac{\delta \mathcal{E}_{h}}{\delta \x_{k}}(\{\x_{i}^{n+1}\}_{i=1}^{N})
\end{split}  
\end{equation}
which can be reformulated as an optimization problem: 
\begin{equation}\label{optimization_problem}
    \argmin_{\{\y_{i}\}_{i=1}^{N}}J(\{\y_{i}\}_{i=1}^{N})  
    = \argmin_{\{\y_{i}\}_{i=1}^{N}}\left[\frac{1}{2\tau}\sum_{i=1}^{N}\|\y_{i} - \x_{i}^{n}\|^{2} + \mathcal{E}_{h}(\{\y_{i}\}_{i=1}^{N}) \right]  
\end{equation}
Hence, we can obtain the particle trajectories by solving the above optimization problem instead of directly solving the ODE system \eqref{ODE_sys}, which may be numerically more expensive.

The Keller–Segel–Patlak system in different dimensions satisfies the general energy dissipation law with the energy functional \eqref{energy} specified by $V = 0$ and
\begin{equation}\label{KS_kernel}
W(x) =
\begin{cases}
-2 \log |x|, & d = 1, \\
\displaystyle \frac{-1}{2\pi} \log |x|, & d = 2, \\
\displaystyle \frac{|x|^{2-d}}{(2-d)\,\omega_{d}}, & d \geq 3,
\end{cases}
\end{equation}
where the interaction kernel is scaled by the chemotactic sensitivity parameter $\chi > 0$, i.e., $W \mapsto \chi W$. For $d = 1$, we consider the modified model studied in \cite{calvez2006modified}, and for d $\geq 3, \omega_{d} = \frac{2\pi^{d/2}}{\Gamma(\frac{d}{2})}$ denotes the surface area of the unit sphere in $\mathbb{R}^{d}$.

The evolution of the system can be interpreted as a competition between two effects: the nonlocal interaction, which promotes aggregation, and the diffusion, which counteracts this by spreading the density. This balance determines whether solutions remain regular or develop concentration phenomena such as blow-up.

\textit{Remark 2.1}
The initial total mass plays a crucial role in determining the aggregation dynamics of the Keller–Segel model \cite{calvez2012blow, blanchet2006two, corrias2004global}. Variations in the initial mass influence the strength and direction of particle interactions, thereby affecting the density evolution and the system’s tendency toward aggregation or dispersion. This dependence is highly sensitive to the spatial dimension.

The most classical case arises when $d = 2$. If the initial mass is below the critical threshold $M_{c} = 8\pi$ (for $\chi = 1$), the solution exists globally in time \cite{blanchet2006two, horstmann2005boundedness}. In contrast, if the initial mass exceeds $M_{c}$, the solution blows up in finite time \cite{blanchet2006two, calvez2012blow, horstmann2005boundedness}.

In one dimension, aggregation can be effectively controlled by tuning the sensitivity parameter $\chi$, whereas in higher dimensions, the same mass and sensitivity can lead to qualitatively different outcomes, including finite-time blow-up \cite{corrias2004global}.
Due to the singularity of W at the origin, we cannot directly apply \eqref{ODE_sys} to compute particle trajectories. Instead, we introduce a regularization of the kernel gradient $\nabla W$. Since $W(x)$ is radially symmetric, we can write
$\nabla W(x) = \phi(|x|)\, x$,
where $\phi$ is the radial profile function. We regularize $\phi$ using a cut-off function $\tilde{\phi}$, defined by
\begin{equation}\label{cut-off_profile}
\tilde{\phi}(|x|) =
\begin{cases}
\phi(|x|), & |x| > r_c, \\
\phi(r_c), & |x| \le r_c,
\end{cases}
\end{equation}
for a chosen cut-off radius $r_c>0$. Therefore, by solving the optimization problem \eqref{optimization_problem} for the Keller–Segel model with $\nabla W(x)$ replaced by $\nabla \tilde{W}(x) := \tilde{\phi}(x) \x$ at each time step, we can obtain the particle trajectories, which serve as our learning data.

\textit{Remark 2.2} 
In contrast to the work \cite{carrillo2019blob}, which derives the particle scheme from the Wasserstein gradient flow structure, our formulation is based on the energetic variational framework. Although the derivations originate from different variational principles—optimal transport versus energy-dissipation balance—the resulting particle dynamics are mathematically equivalent.

Deriving particle schemes from the energetic variational approach offers several advantages. First, it naturally preserves the gradient flow structure, providing a clear energetic interpretation. Second, this framework allows us to bypass the strong regularity assumptions on the density typically required in PDE formulations, enabling particle-based implementation under weaker conditions. Moreover, particle methods are well-suited for high-dimensional problems due to their mesh-free nature.

However, this approach also presents challenges. Incorporating boundary conditions from the continuum PDE into a particle framework is nontrivial. Additionally, aggregation behavior can lead to particle clustering, which may distort the approximation of interaction energy or entropy-like terms and introduce unphysical artifacts in the computed energy. Finally, similar to the scheme derived via the Wasserstein gradient flow, the assumption that the velocity field is evaluated along particle trajectories (i.e., $\dot{\x}_k(t) = \bu(\x_k(t), t)$) may result in an inaccurate energy dissipation rate, particularly in stiff or singular regimes.    
%
\subsection{Stochastic Particle Methods for 2D Keller Segel Model\label{sec:stochastic}}
From the energy–dissipation identity, the two-dimensional parabolic–elliptic Keller–Segel system is given by \eqref{Keller_Segel_equation}:
\begin{equation*}
\partial_t \rho = \Delta \rho + \chi \nabla \cdot \big( \rho , \nabla (W * \rho) \big).
\end{equation*}
From the perspective of stochastic differential equations (SDEs), the dynamics of the Keller–Segel model can be interpreted as a combination of deterministic drift and stochastic diffusion. The nonlocal interaction term, $\tfrac{\chi}{2} W * \rho$, generates the deterministic drift in the associated SDE. This interaction models attractive forces between particles or organisms, such as chemotactic movement toward regions of higher chemical concentration. On the other hand, the diffusion term arising from the entropy functional $\int \rho \log \rho$ corresponds to the stochastic component of the dynamics. It represents random motion (modeled by Brownian motion) and promotes the spreading of the density. 

The interacting particle system $\{\x_{i, t} \in \R^d\}_{i=1}^{N}$ associated with the $2d$ \newline stochastic Keller–Segel model is governed by the SDEs \cite{fournier2017stochastic,misiats2022global,chavanis2010stochastic,cattiaux20162}
\begin{equation}
\dx_{i, t} = \frac{\chi}{N}\sum_{j\neq i} \nabla W\left(\x_{i, t} - \x_{j, t}\right)\dt + \sqrt{2} \dif \mathbf{B}_{i, t}, \quad  i = 1, \ldots, N,
\end{equation}
where $\{B^i_t\}_{i=1}^N$ are independent two-dimensional Brownian motions and the interaction kernel $\nabla W$ is given by
\begin{equation}
\nabla W(\x) = 
\begin{cases}
 -\frac{\x}{2\pi |\x|^{2}}, \quad \x \neq \zero\\
 0,  \quad \x = \zero.
\end{cases}
\end{equation}
To handle the singular interactions, we introduce a regularized particle system \cite{fournier2017stochastic}. The dynamics for each particle $i = 1, \dots, N$ are
\begin{equation}\label{reg_SDE_sys}
\dx^{\epsilon}_{i,t} = \frac{\chi}{N} \sum_{j=1}^{N} \nabla W^{\epsilon}(\x^{\epsilon}_{i,t} - \x^{\epsilon}_{j,t})\dt + \sqrt{2} \dif \mathbf{B}_{i,t}, 
\end{equation}
where $\{\mathbf{B}^i_t\}_{i=1}^{N}$ are two-dimensional independent two-dimensional Brownian motions, and a regularized $\nabla W^{\epsilon}(\x) = -\frac{\x}{2\pi (|\x|^{2} + \epsilon^{2})}$ and $\epsilon >0$ is a small parameter. We use system \eqref{reg_SDE_sys} to learn the regularized profile function $\phi^{\epsilon}(r) = \frac{1}{2\pi(r^{2} + \epsilon^{2})}$. To better capture the aggregation dynamics, the stochastic component is scaled by a small parameter $\eta$, since the unscaled noise would otherwise dominate the particle interactions. Thus, the stochastic system for each particle $i = 1, \dots, N$ is given by
    \begin{equation}\label{scaled_reg_SDE_sys}
\dx^{\epsilon}_{i,t} = \frac{\chi}{N} \sum_{j=1}^{N} \nabla W^{\epsilon}(\mathbf{x}^{\epsilon}_{i,t} - \x^{\epsilon}_{j,t})\dt + \sqrt{2} \eta\dif \mathbf{B}_{i,t}.
\end{equation}
\section{Learning Framework}\label{sec:learn}
Particle trajectories offer a rich source of information about the dynamics of a system, as they track the evolution of particles over time. These trajectories encode the interaction patterns between particles, which provide the data needed to infer the underlying interaction kernel. 
\subsection{A Unified Framework for Learning Interaction Kernels}
In this section, we propose a unified learning framework for inferring the profile function from either deterministic or stochastic approximation of the Keller-Segel PDE.  This approach formulates the problem as a variational inverse problem, where the objective is to learn the profile function as the kernel’s gradient that best captures the system’s dynamics.  Both models fit the standard form of the following SDE of the systems (for $\x_i \in \R^D$, $i = 1, \cdots, N$),
\[
\dx_i(t) = \big(\f(\x_i) + \frac{1}{N}\sum_{j = 1, j \neq i}^N \phi(|\x_j(t) - \x_i(t)|)(\x_j(t) - \x_i(t))\big)\dt + \sigma\dif B_i(t).
\]
Here the constant $\sigma = 0$ gives the deterministic particle method and $\sigma > 0$ gives the stochastic approximation.  We also introduce some vectorized notations
\[
\X = \begin{bmatrix} \x_1 \\ \vdots \\ \x_N \end{bmatrix}, \quad \f_D(\X) = \begin{bmatrix} \f(\x_1) \\ \vdots \\ \f(\x_N) \end{bmatrix}, \quad \F_\phi(\X) = \begin{bmatrix} \frac{1}{N}\sum_{j = 2}^N \phi(|\x_j - \x_1|)(\x_j - \x_1)\\ \vdots \\ \frac{1}{N}\sum_{j = 1}^{N - 1} \phi(|\x_j - \x_N|)(\x_j - \x_N)\end{bmatrix}.
    \]
Here $\X, \f_D, \F_{\phi} \in \R^{D = Nd}$. Furthermore, we introduce an average $\ell_2$ norm on $\R^D$ as
\[
<\X, \Y>_N = \frac{1}{N}<\x_i, \y_i>_2, \quad \X = \begin{bmatrix} \x_1 \\ \vdots \\ \x_N \end{bmatrix} \quad \text{and} \quad \Y = \begin{bmatrix} \y_1 \\ \vdots \\ \y_N \end{bmatrix}.
\]
Therefore $|\X|_N = \sqrt{<\X, \X>_N}$. Our objective is to recover the regularized profile function $\phi^{\text{reg}}$. This function appears in the regularized interaction force $\nabla \tilde{W}(\x) = \phi^{\text{reg}}(r = |\x|)\, \x$. To approximate the unknown function $\phi^{\text{reg}}$, we represent it as a linear combination of basis functions $\{\psi_k\}_{\eta=1}^{n}$:
\begin{equation}\label{profile_representation}
\phi^{\text{reg}}(r) \approx \varphi(r) = \sum_{\eta=1}^{n} \alpha_{\eta} \psi_{\eta}(r),
\end{equation}
where $\alpha_{\eta}$ are the coefficients to be learned. The choice of basis functions depends on prior knowledge of the system and computational considerations. In this paper, we use clamped B-splines as basis functions, motivated by the assumption that the energy is continuous. However, other basis choices, such as the Fourier basis or piecewise polynomial basis, could be employed.  With the representation in \eqref{profile_representation}, the learning task is formulated as a variational inverse problem using the particle trajectories as training data. Specifically, given continuous-time observations of the particle trajectories, i.e., $\{\x_i(t)\}_{i=1}^N$ for $t \in [0,T]$, we define a continuous loss functional as
\begin{equation}\label{continuous_loss_functional}
\mathcal{L}(\varphi) = \tfrac{1}{2}\E\Big[\int_0^T\big(|\F_{\varphi}(\X_t)|_N^2\dt - 2<\F_{\varphi}(\X_t), \dif\X_t - \f_D(\X_t)\dt>_N\big)\Big].
\end{equation}
In practice, we are only given discrete data, i.e., $\{x_i^m(t_l)\}_{i, l, m = 1}^{N, L, M}$, with $M$ denoting the number of initial data distributions, $L$ the number of discrete time points ($0 = t_1 < \cdots < t_L = T$), and $N$ the number of particles in each distribution. Accordingly, the discrete loss functional is defined as 
\begin{equation}\label{discrete_loss_functional}
\widetilde{\mathcal{L}}(\varphi)  = \frac{1}{2ML} \sum_{m,l}^{M, L-1}\big(|\F_{\varphi}(\X_l^m)|_N^2\Delta t - 2<\F_{\varphi}(\X_l^m), \Delta\X_l^m- \f_D(\X_l^m)\Delta t>_N\big)
\end{equation}
where $\X_l^m = \X^m(t_l)$, $\Delta\X_l^m = \X_{l + 1}^m - \X_l^m$, and $\Delta t = t_{l+1} - t_{l}$ for $l = 0, \cdots, L-1$.  Here $\f(\x_{i}^{m}(t_{l}))$ is given by 
\[
\f_D(\x_{i}^{m}(t_{l})) = 
    - \sum_{j=1}^{N} \frac{\nabla K_{h}(\br^{m}_{i,j,l})}{\sum_{k=1}^{N} K_{h}(\br^{m}_{k,i,l})} - \frac{\sum_{j=1}^{N} \nabla K_{h}(\br^{m}_{i,j,l})}{\sum_{k=1}^{N} K_{h}(\br^{m}_{k,i,l})}, \quad \br^m_{k ,i, l} = \x^m_k(t_l) - \x^m_i(t_l).
\]
for the deterministic case; and $\f_D(\x^{m}_{i}(t_{l})) = \zero$ for the stochastic particle trajectories. By using \eqref{profile_representation}, we can write the discrete loss functional \eqref{discrete_loss_functional} as
\begin{equation}\label{linear_regression_loss_functional}
\widetilde{\mathcal{L}}(\vec{\alpha}) = \tfrac{1}{2}\vec{\alpha}^\top A \vec{\alpha} - \vec{b}^\top \vec{\alpha} + \vec{c},
\end{equation}
where $\vec{\alpha} = [\alpha_{1}, \dots, \alpha_{n}]^\top$ and $A$ is a symmetric $n \times n$ matrix defined by (for $\eta, \eta' = 1, \cdots, n$)
\[
     A(\eta, \eta') = \frac{\Delta t}{MLN} \sum_{m,l,i=1}^{M,L-1,N} \frac{1}{N^{2}}  \Bigg\langle \sum_{j=1}^{N} \psi_{\eta}(r^{m}_{i,j,l}) \,  
    \mathbf{r}^{m}_{i,j,l} \ ,  \sum_{k=1}^{N} \psi_{\eta'}(r^{m}_{i,j,l}) \, \mathbf{r}^{m}_{i,j,l}\Bigg\rangle_{2}, 
\]
where $r^{m}_{i,j,l} = \|\mathbf{r}^{m}_{i,j,l}\|_{l^{2}(\mathbb{R}^{d})}$. The vector $\vec{b}$ is defined by (for $\eta = 1, \cdots, n$)
\[
    \vec{b}(\eta) = \frac{1}{MLN} \sum_{m,l,i=1}^{M,L-1,N} \Bigg\langle \Delta_{l}\x^{m}_{i}  - \f(\x^{m}_{i}(t_{l}))\Delta t, 
    \frac{1}{N} \sum_{j=1}^{N}\psi_{\eta}(r^{m}_{i,j,l}) \, \mathbf{r}^{m}_{i,j,l}\Bigg\rangle_{2}, 
\]
and $\vec{c}$ denotes the part of the loss functional that is independent of the coefficient vector $\vec{\alpha}$. Thus the loss functional can be reformulated as a quadratic form in $\vec{\alpha}$ with $A$ being SPD.  Minimizing \eqref{discrete_loss_functional} over $\vec{\alpha}$ yields a linear system, i.e. $A \vec{\alpha} = \vec{b}$.  Once the optimal coefficient vector $\vec{\alpha}$ is obtained, we reconstruct $\hat\phi(|x|)$ and use it to simulate Keller-Segel particle trajectories. The learned profile function is validated by comparing the generated trajectories with the original data.  We summarize the above discussion into an algorithm:
\begin{algorithm}[H]
\caption{Learning the Regularized Profile Function $\phi^{\text{reg}}(|x|)$}
\textbf{Input:} Particle trajectories $\{\x_i^m(t_l)\}_{i=1}^N, \; l=1,\dots,L, \; m=1,\dots,M$; time step $\Delta t$; basis functions $\{\psi_\eta\}_{\eta=1}^n$. \\  
\textbf{Output:} Estimated profile function $\hat\phi(|x|)$.

\textbf{Step 1:} Represent the test function as 
$\varphi = \sum_{\eta=1}^n \alpha_\eta \psi_\eta(|x|)$. \\ 

\textbf{Step 2:} Construct the discrete loss functional $\widetilde{\mathcal{L}}(\vec{\alpha})$
based on \eqref{discrete_loss_functional}.  \\ 

\textbf{Step 3:} Assemble the matrix $A \in \mathbb{R}^{n \times n}$ and vector 
$\vec{b} \in \R^n$ as defined in \eqref{linear_regression_loss_functional}.  \\ 

\textbf{Step 4:} Solve the linear system $A\vec{\alpha} = \vec{b}$ to obtain $\vec{\hat\alpha} = [\hat\alpha_1,\dots,\hat\alpha_n]^\top$. \\

\textbf{Step 5:} Reconstruct the profile function $\hat\phi(|x|) = \sum_{\eta=1}^n \hat\alpha_\eta \psi_\eta(|x|)$. \\

\textbf{Step 6:} Simulate Keller–Segel particle trajectories using $\nabla \hat{W}(x) = \hat\phi(|x|)\,x$. \\

\textbf{Step 7:} Validate by comparing simulated trajectories with the observed data. 

\end{algorithm}
%
\subsection{Performance Measures}\label{subsec:3A}
To quantify the accuracy of the learned interaction profile and the reconstructed particle trajectories with the regularized function  $\phi^{\text{reg}}$ and the original particle trajectory data, we define two complementary error measures: the trajectory error and the radial profile error. First, to quantify the discrepancy between the original data and the reconstructed trajectories, we define the numerical trajectory error as follows. Let $\X$ be the system trajectory defined as before, then denote the state variable of the system at time $t$ for the $m^{th}$ initial condition. The trajectory error for each initial distribution is given by
\begin{equation}
\text{Err}_{\X}^m  = \left(\int_{0}^{T} \|\X^{m}_{t} - \hat{\X}^{m}_{t}\|^{2}_{N}dt\right)^{1/2}, 
\end{equation}
which we approximate numerically as 
\begin{equation}
\text{Err}_{\X}^m \approx \left(\frac{1}{L}\sum_{l=0}^{L-1}\|\X^{m}_{t_{l+1}} - \hat{\X}^{m}_{t_{l+1}}\|^{2}_{N}(t_{l+1} - t_{l})\right)^{1/2},  
\end{equation}
for the observation points $\{t_{l}\}_{l=0}^{L}$.  Here, $\X^{m}_{t} $ represents the original particle trajectories for the $m^{th}$ initial condition at time $t$ and $\hat{\X}^{m}_{t}$ denotes the trajectories reconstructed by substituting the learned profile function $\hat\phi$ into \eqref{particle_evolution_deter} and \eqref{reg_SDE_sys}. Assuming an equidistant time step  $\Delta t = t_{l+1} - t_{l}, \ l = 0, \cdots, L-1$, we compute the overall error statistics across all initial conditions as
\begin{equation}\label{trajectory_error}
\begin{split}
    \text{Err}_{\X}  = \frac{1}{M}\sum_{m=1}^{M} \left[\sum_{l=0}^{L-1}\|\X^{m}_{t_{l+1}} - \hat{\X}^{m}_{t_{l+1}}\|^{2}_{N}\Delta t \right]^{1/2},
\end{split}
\end{equation}
where $M$ is the total number of the given initial data distributions. Therefore, the relative trajectory data error can be defined as 
\begin{equation}\label{rel_trajectory_error}
\begin{split}
\text{Err}^{\text{rel}}_{\X}  = \dfrac{\frac{1}{M}\sum_{m=1}^{M} \left[\sum_{l=0}^{L-1}|\X^{m}_{t_{l+1}} - \hat{\X}^{m}_{t_{l+1}}|^{2}_{N}\Delta t \right]^{1/2}}{\frac{1}{M}\sum_{m=1}^{M} \left[\sum_{l=0}^{L-1}|\X^{m}_{t_{l+1}}|^{2}_{N}\Delta t \right]^{1/2}}.
\end{split}
\end{equation}
This measure provides an aggregate assessment of the trajectory reconstruction error, reflecting the deviation of the reconstructed trajectories from the original data.

The $ L^2(\rho)$ error of comparing $\phi^{\text{reg}}(r)r$ is given by
\begin{equation}\label{data_error}
\|\phi^{\text{reg}}(\cdot)\cdot - \hat\phi(\cdot)\cdot \|^{2}_{L^{2}(\rho)} = \int_{r = a}^{b} |\phi^{\text{reg}}(r) - \hat\phi(r)|^{2} \ r^{2} \, \dif\rho(r). 
\end{equation}
Here the data distribution $\rho$ is given by
\[
\rho(r) = \E\left[\frac{2}{N(N - 1)}\sum_{1 \le i < j \le N}\frac{1}{T}\int_0^T\delta_{r_{i, j}(t)}(r) \, \dif r\right], \quad r_{i, j}(t) = \norm{\x_i(t) - \x_j(t)}.
\]
Thus, the relative discrete radial profile error can be defined as 
\begin{equation}\label{rel_profile_error}
    Err_{\phi}^{\text{rel}} = \frac{\|\phi^{\text{reg}}(\cdot)\cdot - \hat\phi(\cdot)\cdot \|^{2}_{L^{2}(\rho)}}{\|\phi^{\text{reg}}(\cdot)\cdot\|^{2}_{L^{2}(\rho)}}
\end{equation}
\subsection{Adaptive Learning}\label{sec:adaptive_learn}
\hspace*{1em} In many applications, including trajectory reconstruction and kernel learning, the distribution of data points often exhibits a nonuniform nature. This necessitates an adaptive approach to selecting knot points for the B-spline basis functions to achieve an efficient and accurate representation of the underlying function. The use of nonuniform knots allows the spline to focus more on regions of the domain where data points are densely concentrated, while reducing resolution in regions with sparse data.

In standard B-spline formulations, knots are typically chosen to be uniformly spaced. However, for complex data distributions, nonuniform knot placement can lead to better accuracy and efficiency—particularly when certain regions exhibit significant variation or require higher resolution. Building on these observations, we introduce an adaptive procedure that iteratively refines knot placement based on local error estimates. We first define minimum and maximum pairwise distances for the given particle trajectory data as 
\begin{equation}\label{initial_partition}
    a = \underset{l,i,j,m}{\min}r^{m}_{i,j,l}, \quad  b = \underset{l,i,j,m}{\max}r^{m}_{i,j,l},  
\end{equation}
where $r^{m}_{i,j,l} = \|x^{m}_{i}(t_{l}) - x^{m}_{j}(t_{l})\|_{l^{2}(\mathbb{R}^{d})}$. The algorithm to refine the knot placement is as follows: 
\begin{algorithm}[H]
\textbf{Step 1}: Define the initial partition and indicator set:
\begin{equation}
\begin{cases}
\mathcal{P}^{0} = \{a, b\}, \\
\mathcal{I}^{0} = \{1\}.
\end{cases}
\end{equation}

For each iteration counter $j$, set:
\begin{equation}
\begin{cases}
\mathcal{I}^{j} = \{B_1, \dots, B_K\}, \\
\mathcal{P}^{j} = \{x_0^{j}, \dots, x_K^{j}\},
\end{cases}
\end{equation}
where $B_k$ is a boolean indicating whether $[x_{k-1}, x_k]$ should be divided. Build new partitions and solve intermediate problems:
\begin{equation}
\begin{cases}
\mathcal{P}^{j-1} = \{x^{j-1}_{0}, \dots, x^{j-1}_{K^{j-1}}\} \Rightarrow \text{solve for } \hat{\phi}_{j-1}, \\ 
\tilde{\mathcal{P}}^{j} = \{\tilde{x}^{j}_{0}, \dots, \tilde{x}^{j}_{2K^{j-1}}\} \Rightarrow \text{solve for } \hat{\phi}_{j},
\end{cases}
\end{equation}
with:
\begin{equation}
\tilde{x}^{j}_{2k} = x^{j-1}_{k}, \quad 
\tilde{x}^{j}_{2k-1} = \frac{x^{j-1}_{k-1} + x^{j-1}_{k}}{2}.
\end{equation}

\textbf{Step 2}: Compare errors on each subinterval:
\begin{equation}
E_{rel}  = \frac{\sqrt{\int_{x^{j-1}_{k-1}}^{x^{j-1}_{k}}|\hat{\phi}_{j-1}(r) - \hat{\phi}_{j}(r)|^{2}d\rho(r)}}
{\sqrt{\int_{x^{j-1}_{k-1}}^{x^{j-1}_{k}}|\hat{\phi}_{j-1}(r)|^{2}d\rho(r)}} \approx \frac{\left|\hat{\phi}_{j-1}(\tilde{x}^{j}_{2k-1}) - \hat{\phi}_{j}(\tilde{x}^{j}_{2k-1})\right|}
{\left|\hat{\phi}_{j-1}(\tilde{x}^{j}_{2k-1})\right|}.
\end{equation}

Update indicators:
\[
\begin{cases}
E_{rel} \leq \text{tol} \Rightarrow B^{j-1}_{k} = 0, \\
E_{rel} > \text{tol} \Rightarrow B^{j-1}_{k} = 1.
\end{cases}
\]

Termination:
\[
\begin{cases}
\text{If all } E_{rel} \leq \text{tol} \text{ or } j > \text{maxIter}: \text{ stop and return } \mathcal{P}^{j}, \\
\text{Else: construct } \mathcal{P}^{j} \text{ from } \mathcal{P}^{j-1} \text{ and } \mathcal{I}^{j-1}.
\end{cases}
\]
\caption{Adaptive Learning of Partition Based on Relative $L^2_{r^{2}}(\rho)$ Error}
\end{algorithm}
The concept of adaptive knot selection involves placing more knots in regions where the function or data exhibits greater complexity, and fewer knots where the function is smoother or the data is sparse. This process leverages the local behavior of the function or data, allowing for finer resolution where necessary. Drawing inspiration from error estimation techniques in the finite element method, we iteratively refine the knot placement based on local error estimates, ensuring that the refinement is adaptively tailored to the kernel learning process.
\section{Examples}\label{sec:examples}
In this section, we present numerical results with performance measures for learning the interaction profiles functions. The objective is to assess the accuracy of the learned profile function in capturing the underlying dynamics and the effectiveness of the reconstructed trajectories in approximating true particle motions. We consider various regularization techniques and parameter choices to ensure stability and generalization.  In all the deterministic examples given below, our mollifiers $K_{h}$ are Gaussian functions  
\begin{equation*}
    K_{h}(x) = \frac{1}{(2\pi h^{2})^{d/2}} e^{-|x|^{2}/2h^{2}}.
\end{equation*}
The initial distribution of particle positions for each dataset is uniformly distributed in the unit cube $[0,1]^d$. The partition number of the interval $[a, b]$, as defined in \eqref{initial_partition}, is set to $P = 400$. To approximate the regularized profile functions $\tilde{\phi}$, we use cubic B-splines as our basis functions.  For the adaptive learning, we set the tolerance for the relative error to $\text{tol} = 0.01$ and the maximum number of iterations to $\text{maxIter} = 6$. For both the regularized particle trajectories and the reconstructed particle trajectories using the learned profile $\phi^{K}$, the time step $\tau$ is set to $1 \times 10^{-4}$.
  
We list all the common parameter values used throughout the numerical tests in Table~\ref{tab:common_parameters}.
\begin{table}[h] 
    \centering
    \captionsetup{width=\linewidth}
    \renewcommand{\arraystretch}{1.0}  
    \setlength{\tabcolsep}{8pt}       
    \begin{tabular}{l c c}
        \toprule
        \textbf{Parameter} & \textbf{Symbol} & \textbf{Value} \\  
        \midrule
        Final time & $T$ & 0.2 \\  
        Observation time step & $\Delta t$ & 0.01 \\  
        Time step for particle trajectories & $\tau$ & $1 \times 10^{-4}$ \\ 
        Number of particles & $N$ & 50 \\  
        Number of initial data samples & $M$ & 500 \\  
        Partition number of the interval & $P$ & 400 \\  
        Bandwidth for $K_{h}$ & $h$ & 0.01 \\  
        Tolerance for relative error & $\text{tol}$ & 0.01 \\  
        Diffusion coefficient & $\eta$ & 0.01 \\
        Regularization parameter for $\nabla W^{\epsilon}$ & $\epsilon$ & 0.01 \\ 
        Maximum iterations \\ (adaptive knot points) & $\text{maxIter}$ & 6 \\  
        \bottomrule
    \end{tabular}
    \caption{Common parameter values.}
    \label{tab:common_parameters}
\end{table}
\subsection{The One-Dimensional Case}
In the one-dimensional case, we consider the modified Keller–Segel model studied in \cite{calvez2006modified}. We choose the parameters $\chi = 0.35$, $0.55$, and $0.75$ to assess the accuracy of the learned profile function using uniform knot points. For adaptive knot points, we set $\chi = 0.55$. The truncation parameter $r_c$ is fixed at $0.01$ for both the uniform and adaptive knot cases.

\begin{figure}[ht]
    \centering 
    \captionsetup{width=\linewidth} 
    \begin{subfigure}{0.35\linewidth}
        \centering
        \includegraphics[width=\linewidth]{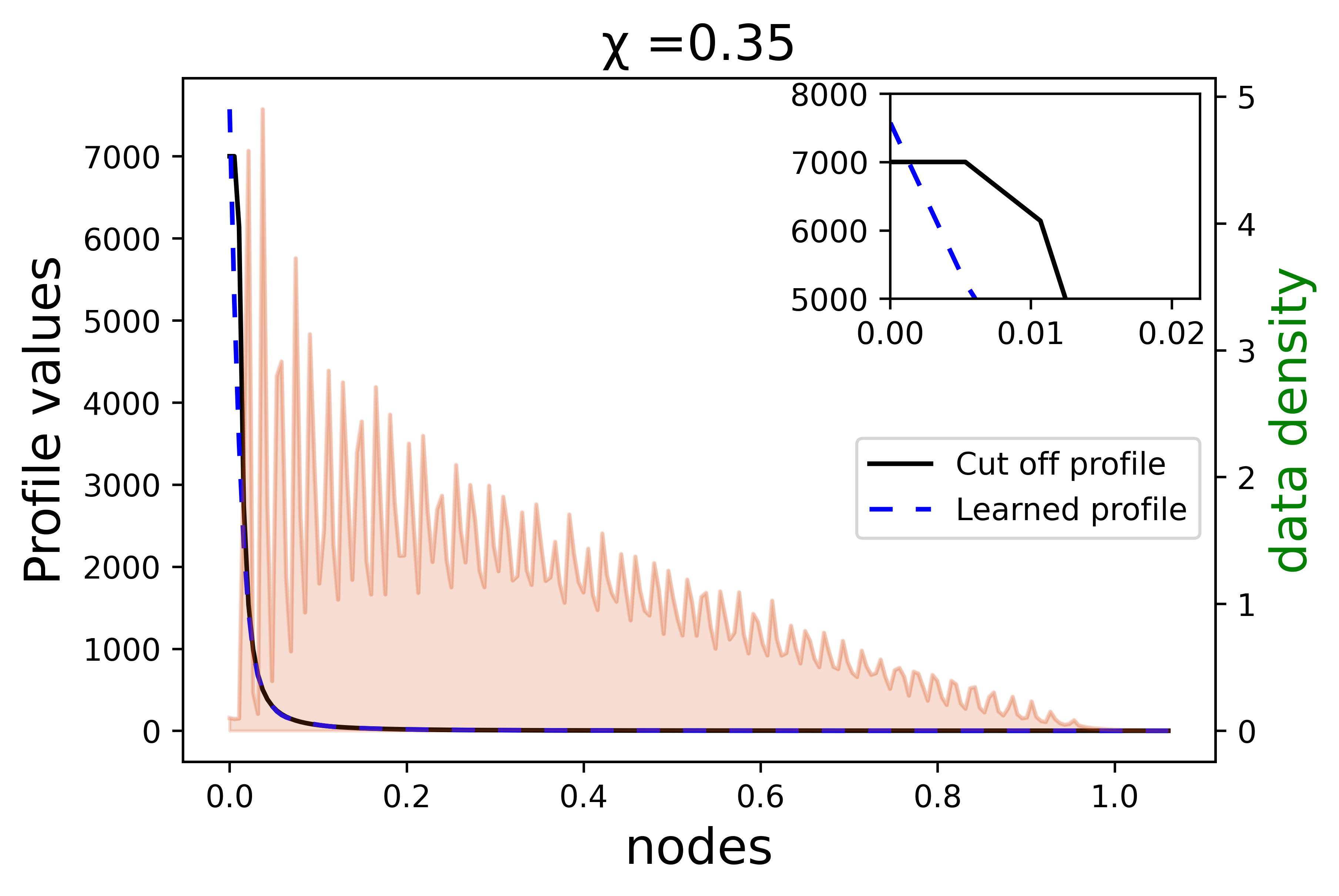}
        \caption{}
    \end{subfigure}
    \begin{subfigure}{0.6\linewidth}
        \centering
        \includegraphics[width=\linewidth]{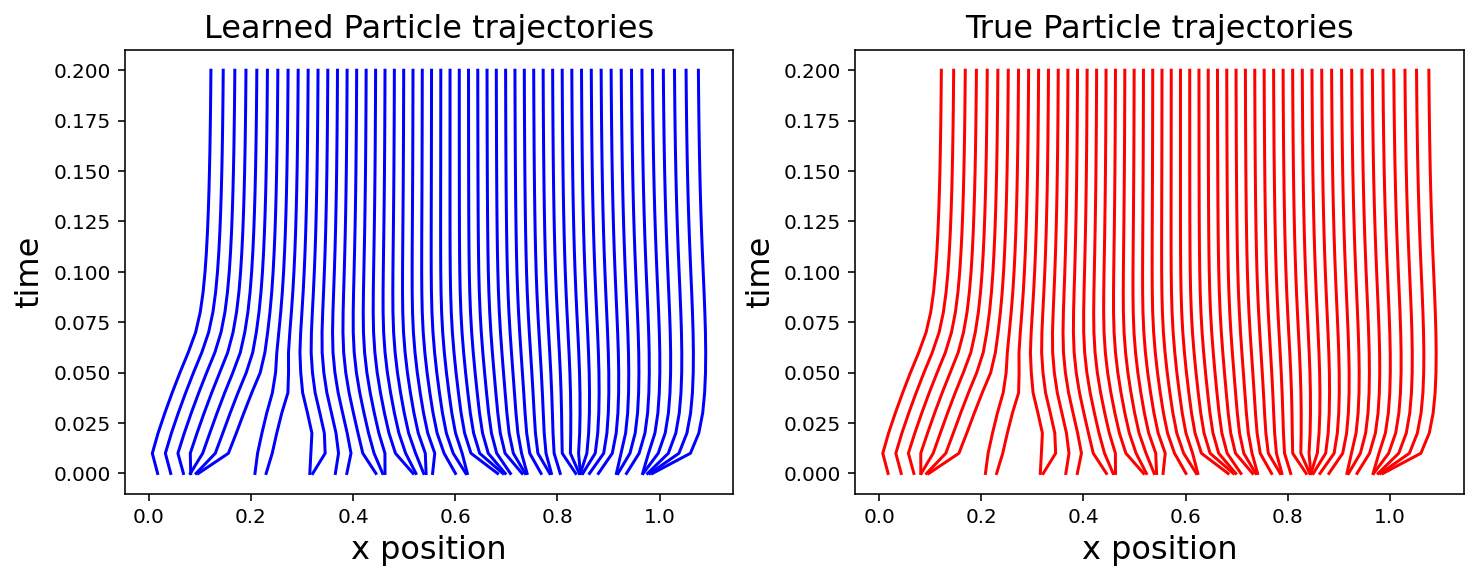}
        \caption{}
    \end{subfigure}  
    \begin{subfigure}{0.35\linewidth}
        \centering
        \includegraphics[width=\linewidth]{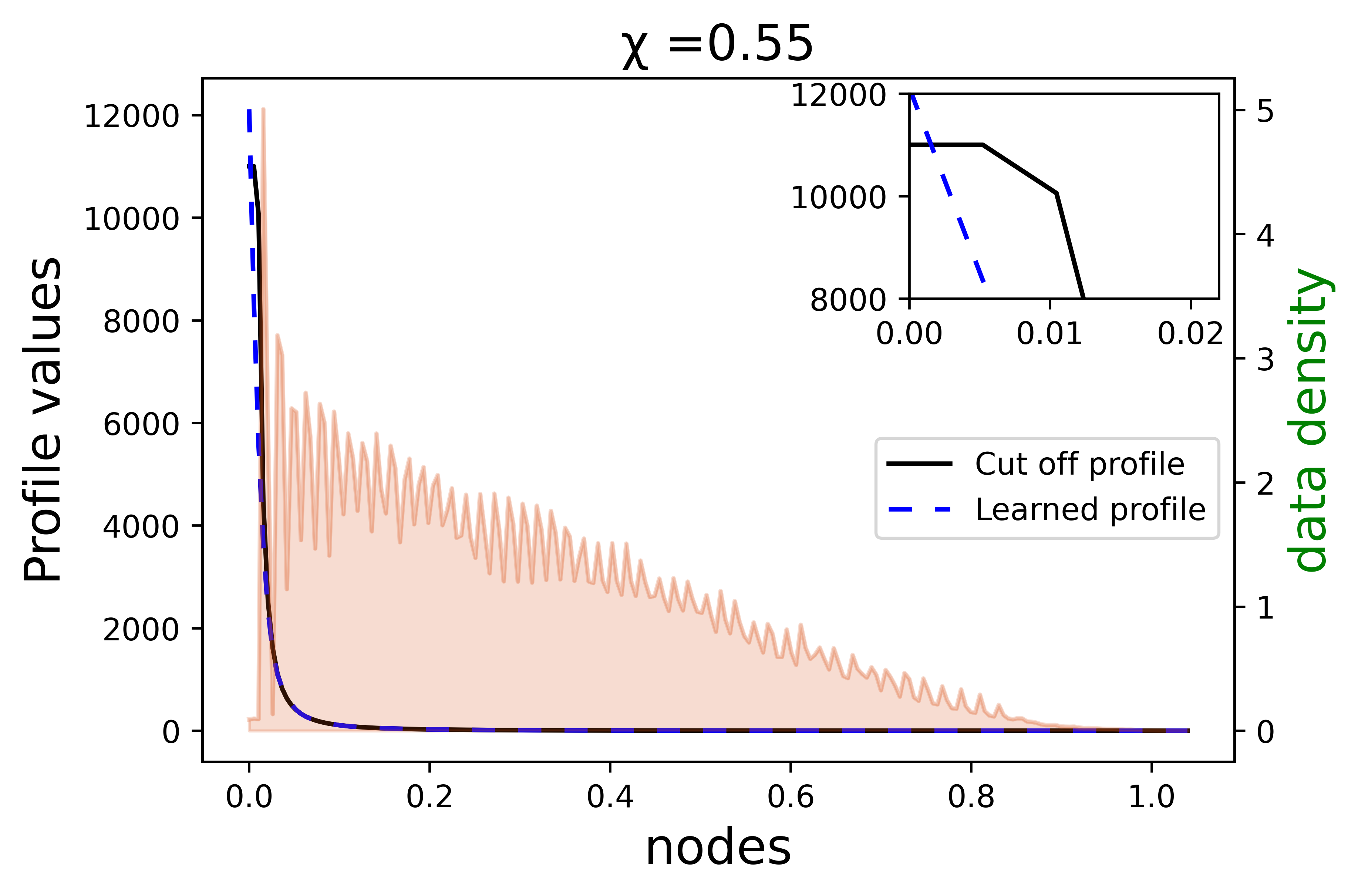}
        \caption{}
    \end{subfigure}
    \begin{subfigure}{0.6\linewidth}
        \centering
        \includegraphics[width=\linewidth]{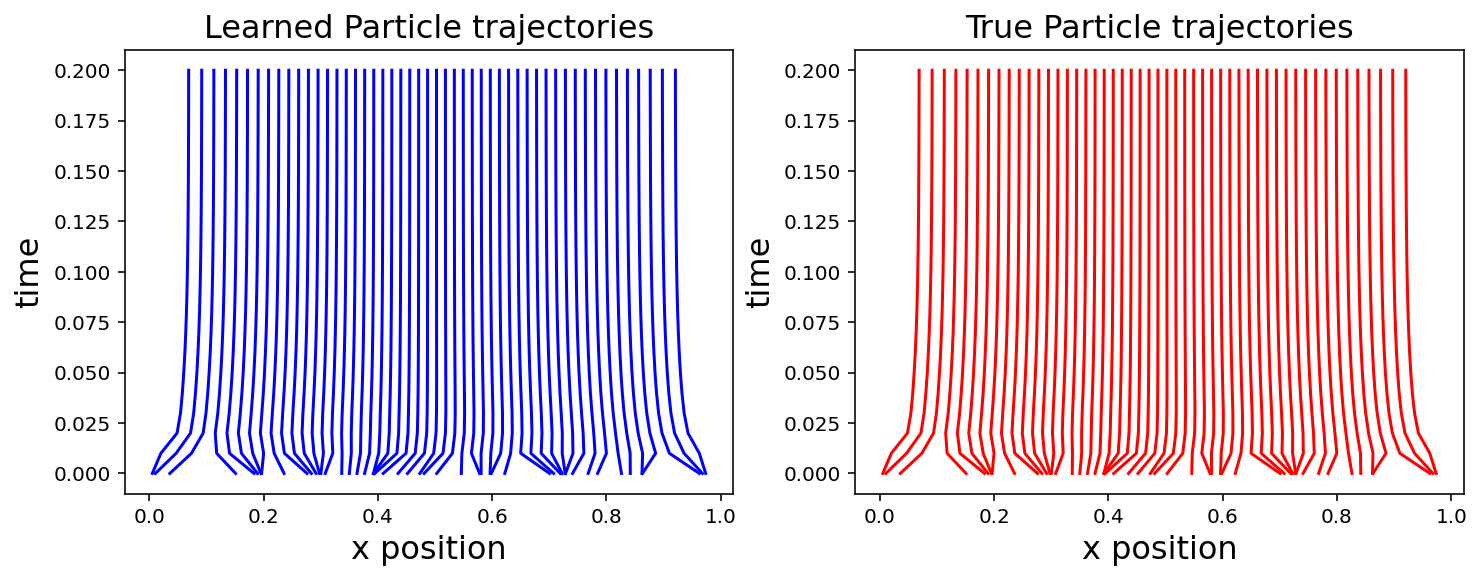}
        \caption{}
    \end{subfigure} 
    
    \begin{subfigure}{0.35\linewidth}
        \centering
        \includegraphics[width=\linewidth]{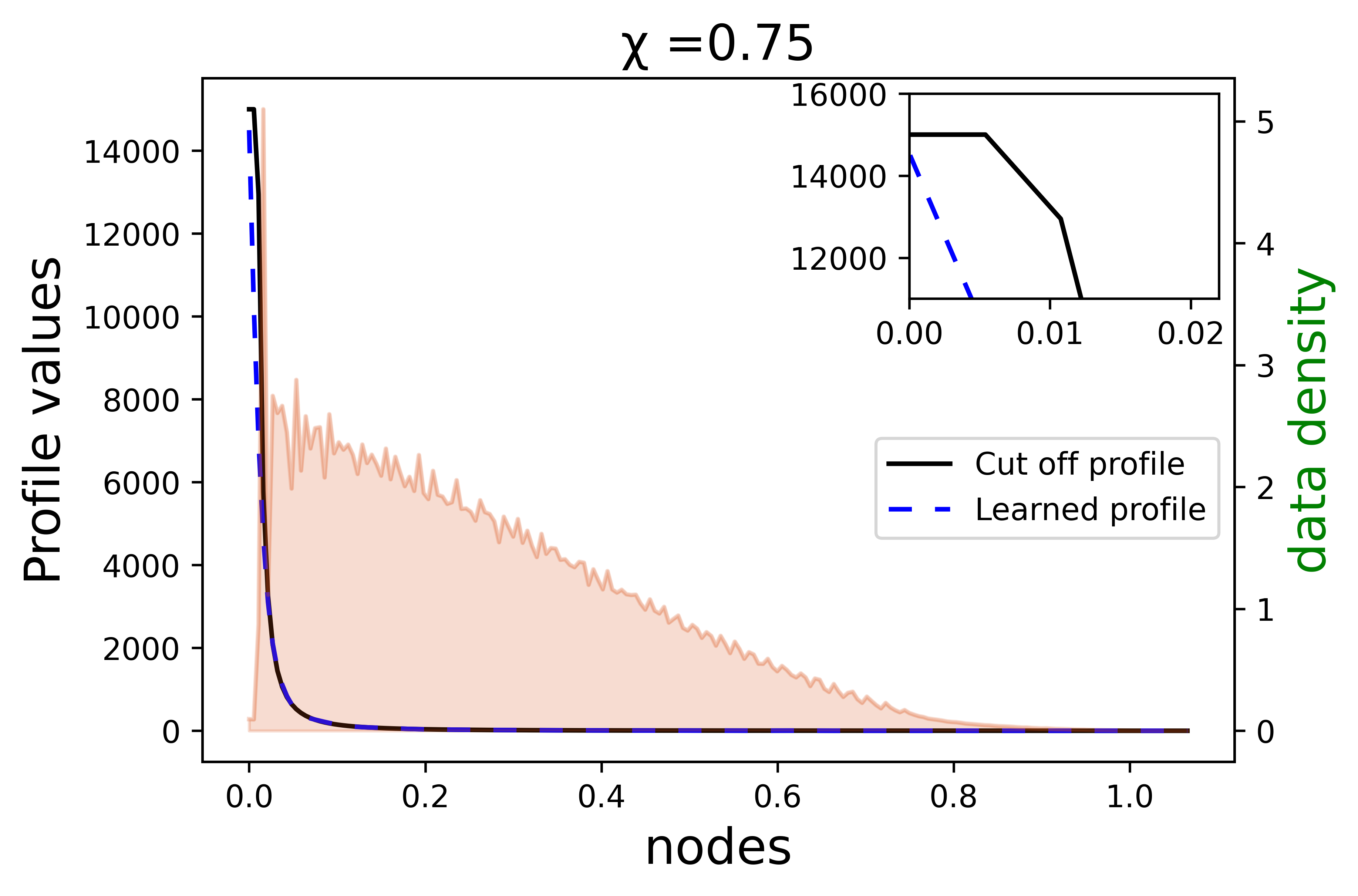}
        \caption{}
    \end{subfigure}  
     \begin{subfigure}{0.6\linewidth}
        \centering
        \includegraphics[width=\linewidth]{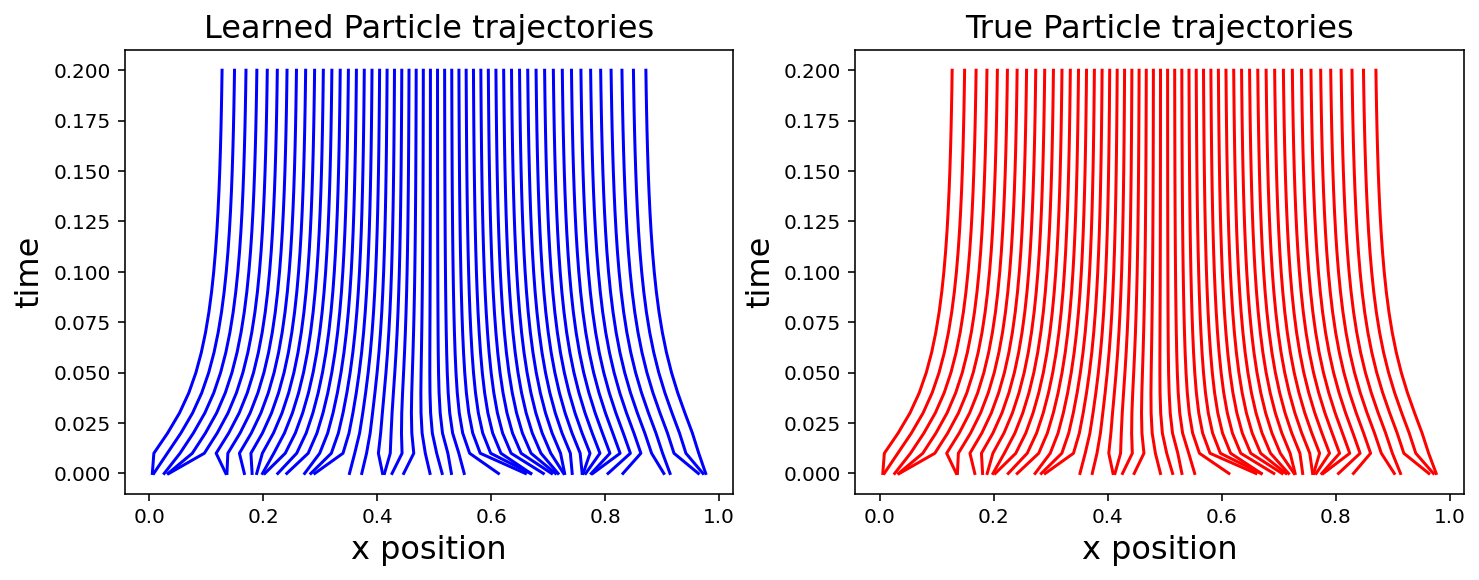}
        \caption{}
    \end{subfigure}
    \caption{\label{fig:1D_kernel_uniform_knot} Comparison of the learned profile functions with the regularized profile functions in the one-dimensional modified Keller–Segel model, using 30 uniform knot points and a truncation parameter of $r_c = 0.01$. The top, middle, and bottom subfigures correspond to the parameters $\chi = 0.35$, $0.55$, and $0.75$, respectively. The corresponding particle trajectories—learned (left) and true (right)—are shown for each $\chi$ value.}
\end{figure}

\begin{figure}[h]
    \centering 
    \captionsetup{width=\linewidth}
    \begin{subfigure}{0.45\linewidth}
        \centering
        \includegraphics[width=\linewidth]{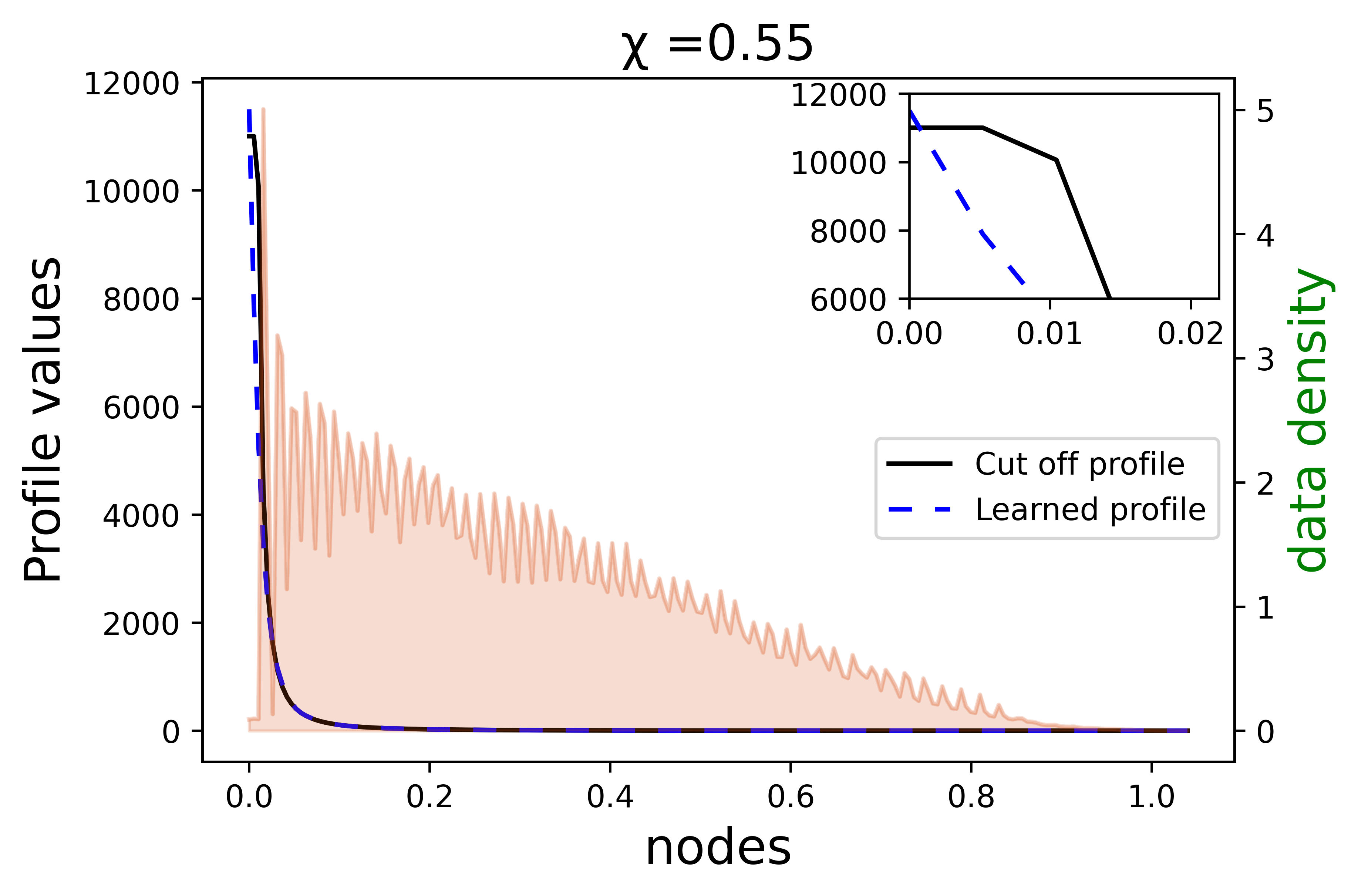}
        \caption{}
    \end{subfigure}
    \hfill 
    \begin{subfigure}{0.45\linewidth}
        \centering
        \includegraphics[width=\linewidth]{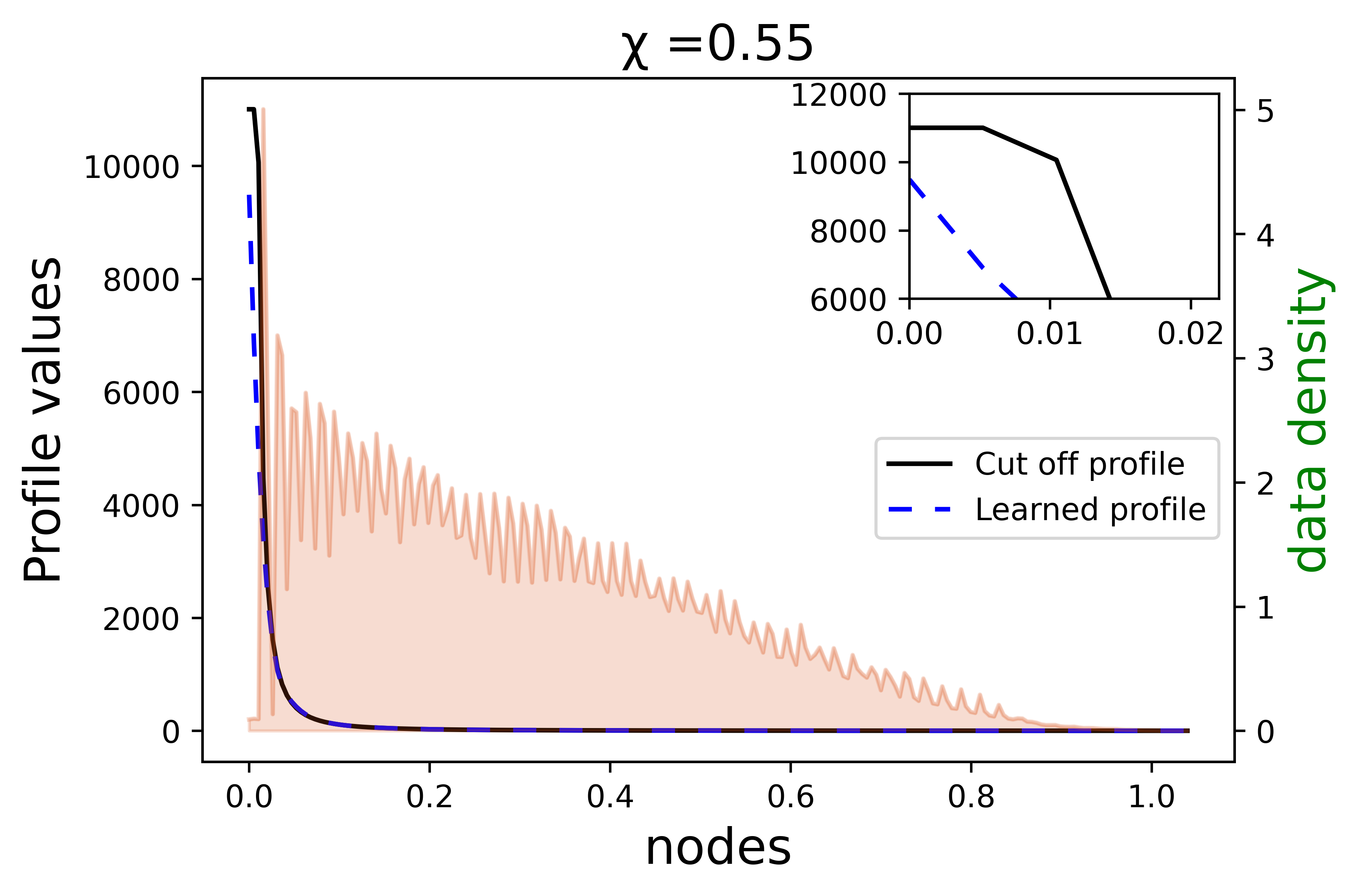}
        \caption{}
    \end{subfigure}
\caption{Comparison of the learned profile functions with the regularized profile functions in the one-dimensional modified Keller–Segel model, using 25 adaptive (left) and uniform (right) knot points, with $r_c = 0.01$ and $\chi = 0.55$. Insets highlight the profile behavior close to the origin.}
\label{fig:1D_kernel_adaptive_knot}
\end{figure}

Figure~\ref{fig:1D_kernel_uniform_knot} compares the learned profile function $\phi^K$ with the regularized profile function $\tilde{\phi}$, both obtained using 30 uniform knot points, and the corresponding reconstructed particle trajectories. In the left subfigure, the background plot represents the data density of pairwise particle distances, which exhibits a decaying trend as the pairwise distance increases. The learned kernel closely matches the regularized profile near the origin, where the data density is higher. In the right subfigure, the reconstructed particle trajectories obtained using the learned profile functions (shown in the left subfigure) are compared with the regularized trajectories. The aggregation phenomenon observed at the beginning, followed by the cessation of particle motion, is well captured by the reconstructed trajectories using the learned profile $\phi^{K}$.  
 \\ 
\indent Figure~\ref{fig:1D_kernel_adaptive_knot} compares the learned profile functions obtained using uniform and adaptive knot points. The number of knots is set to 25 for both cases. The results show that the profile obtained with adaptive knots matches the regularized function more accurately than the one with uniform knots. This suggests that adaptive knot points are better suited for capturing the underlying profile, likely due to their ability to concentrate more on regions with higher variability, resulting in a more accurate approximation of the interaction kernel.

\indent Table~\ref{tab:1D} reports the numerical relative errors in the reconstructed trajectories, $Err_{\text{traj}}^{\text{rel}}$ (defined in \ref{rel_trajectory_error}), and in the learned profiles, $Err_{\phi}^{\text{rel}}$ (defined in \ref{rel_profile_error}), for the one-dimensional case with different values of the sensitivity parameter $\chi$. We observe that the trajectory errors remain small (on the order of $10^{-4}$) and decrease slightly as $\chi$ increases. In contrast, the profile errors grow with larger $\chi$, ranging from $0.096$ at $\chi=0.35$ to $0.207$ at $\chi=0.75$. This indicates that while trajectory reconstruction is robust across parameter values, learning the profile becomes more challenging due to the increasing singularity of the profile function.
\begin{table}[h!]
  \centering 
  \captionsetup{width=\linewidth}  
  \renewcommand{\arraystretch}{1.3} 
  \setlength{\tabcolsep}{6pt} 
  \begin{minipage}{0.95 \columnwidth}
  \centering
  \begin{tabular}{ |c | c c c |}
    \hline
    $\chi$ & $0.35$ & $0.55$ & $0.75$ \\ \hline
    $Err_{\text{traj}}^{\text{rel}}$ & $3.94\mathrm{e}{-4}$ & $3.55\mathrm{e}{-4}$ & $3.35\mathrm{e}{-4}$ \\ \hline 
    $Err_{\phi}^{\text{rel}}$ & $0.096$ & $0.116$ & $0.207$ \\ \hline
  \end{tabular}
  \caption{Performance measures ($1D$ case) using $30$ uniform knots and $r_{c} = 0.01$.}
  \label{tab:1D}
  \end{minipage}
\end{table} 
\subsection{The Two-Dimensional Case}
\begin{figure}[h] 
    \centering 
    \captionsetup{width=\linewidth}
    \begin{subfigure}{0.37\linewidth}
        \centering
        \includegraphics[width=\linewidth]{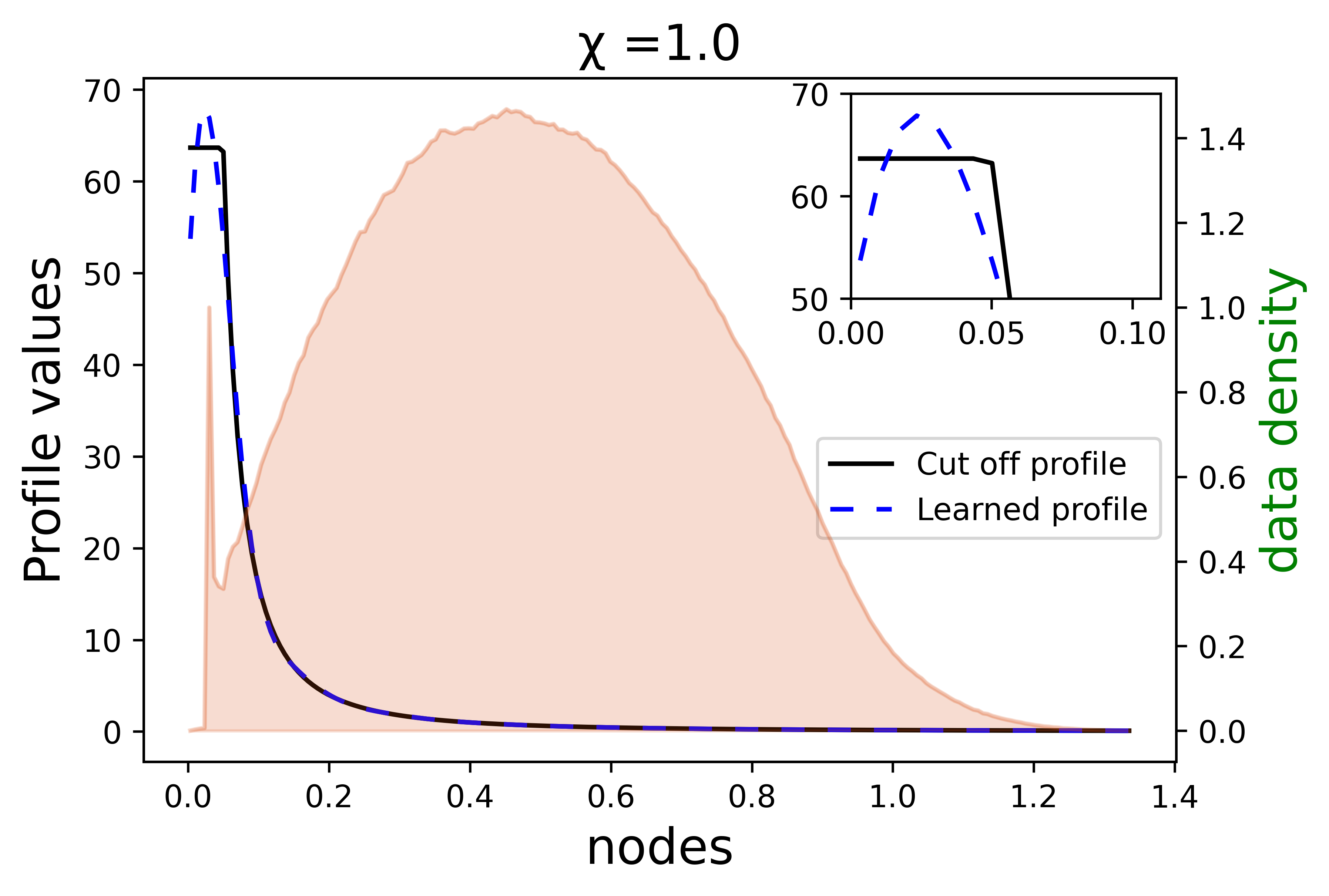}
        \caption{}
    \end{subfigure}
    \begin{subfigure}{0.6 \linewidth}
        \centering
        \includegraphics[width=\linewidth]{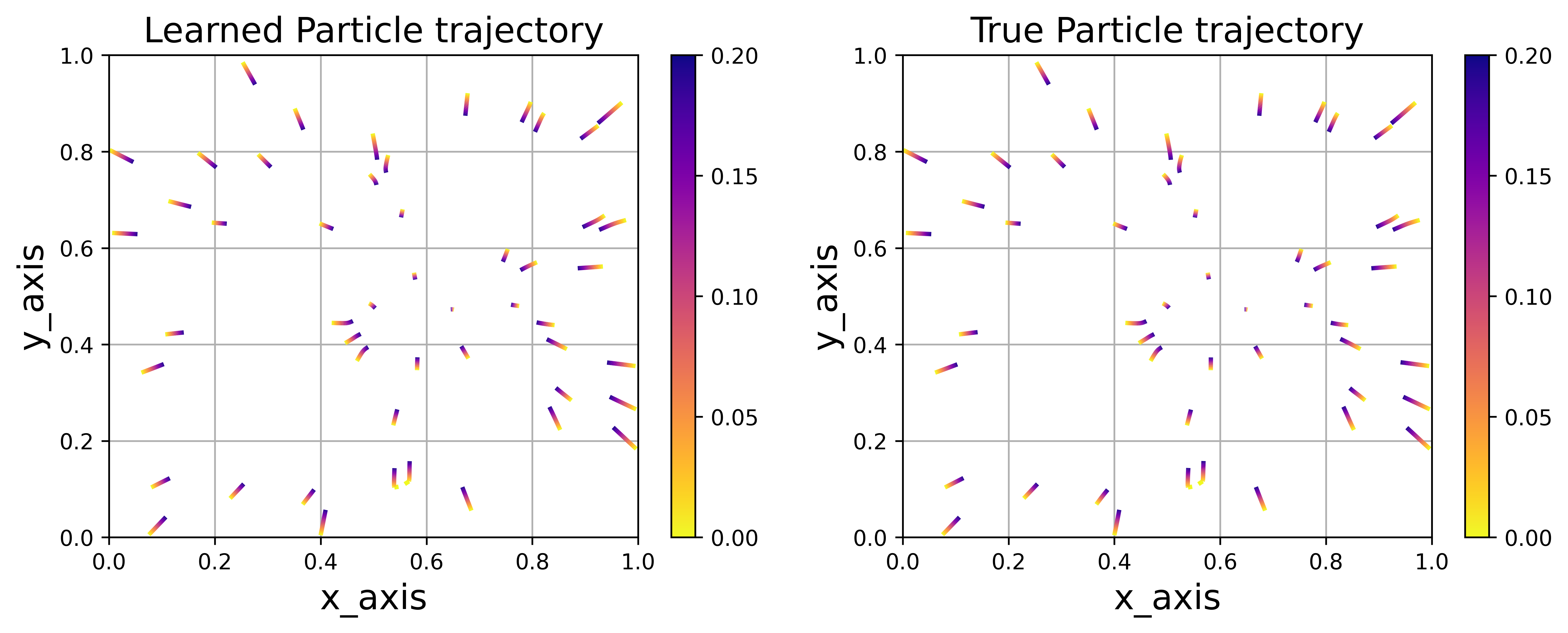}
        \caption{}
    \end{subfigure}
    \begin{subfigure}{0.37\linewidth}
        \centering
        \includegraphics[width=\linewidth]{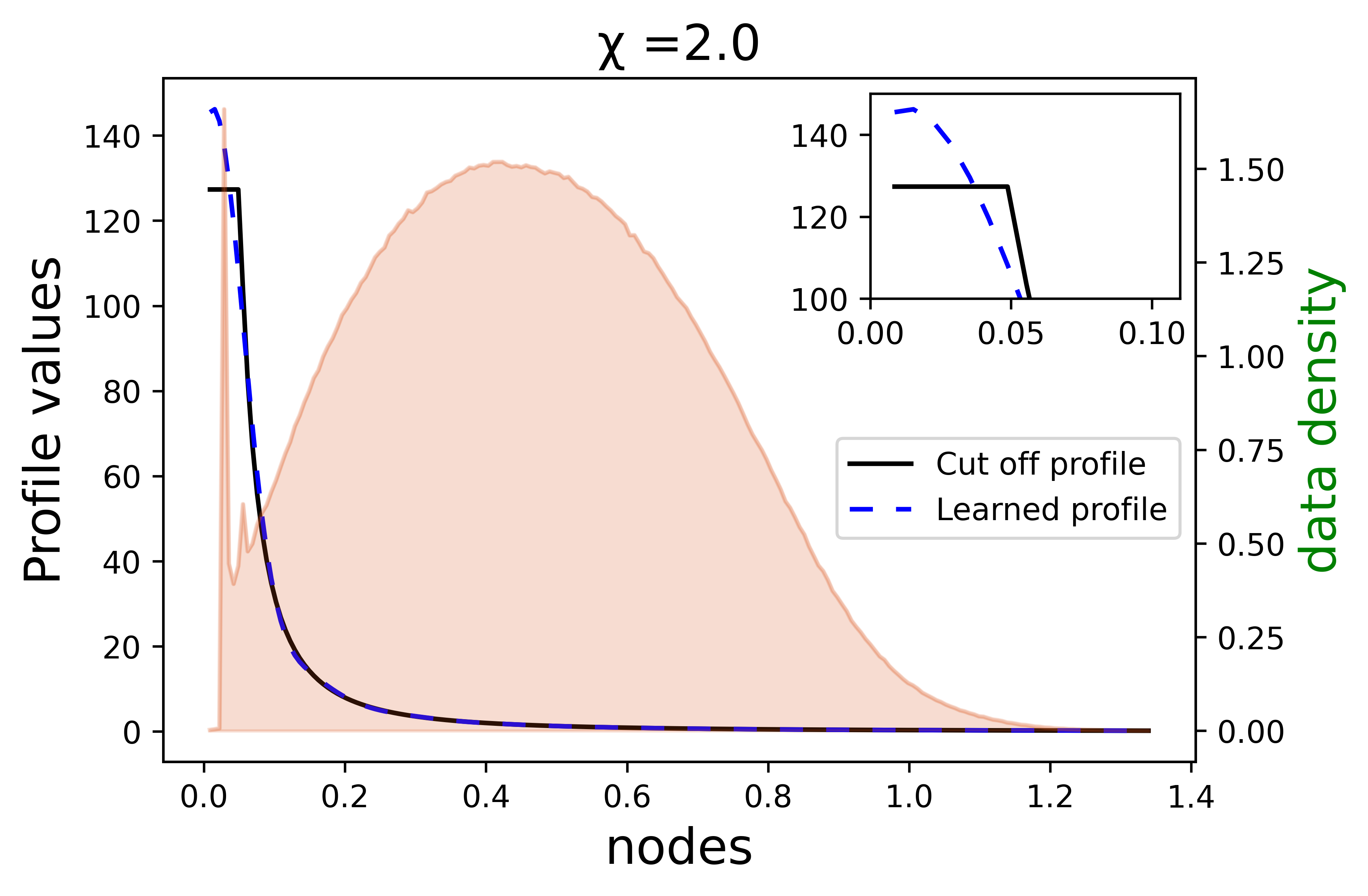}
        \caption{}
    \end{subfigure}
    \begin{subfigure}{0.6\linewidth}
        \centering
        \includegraphics[width=\linewidth]{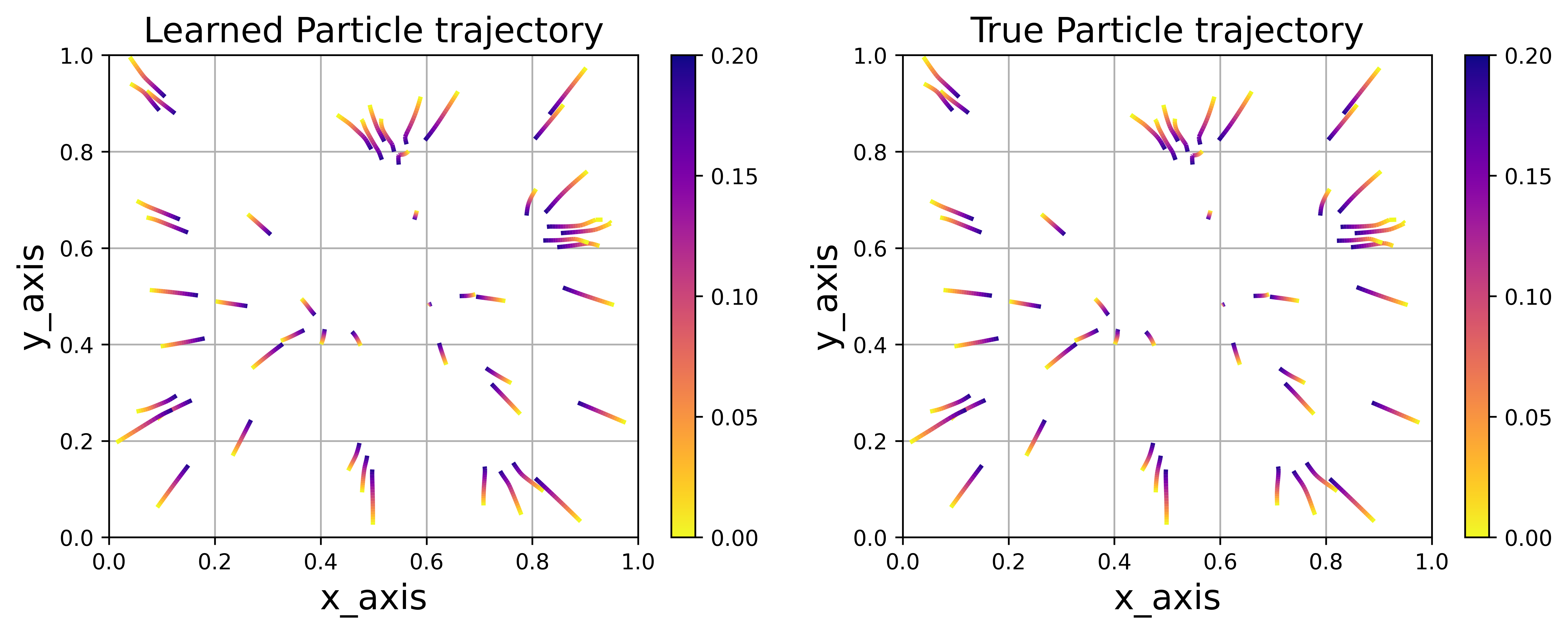}
        \caption{}
    \end{subfigure}
    \begin{subfigure}{0.37\linewidth}
        \centering
        \includegraphics[width=\linewidth]{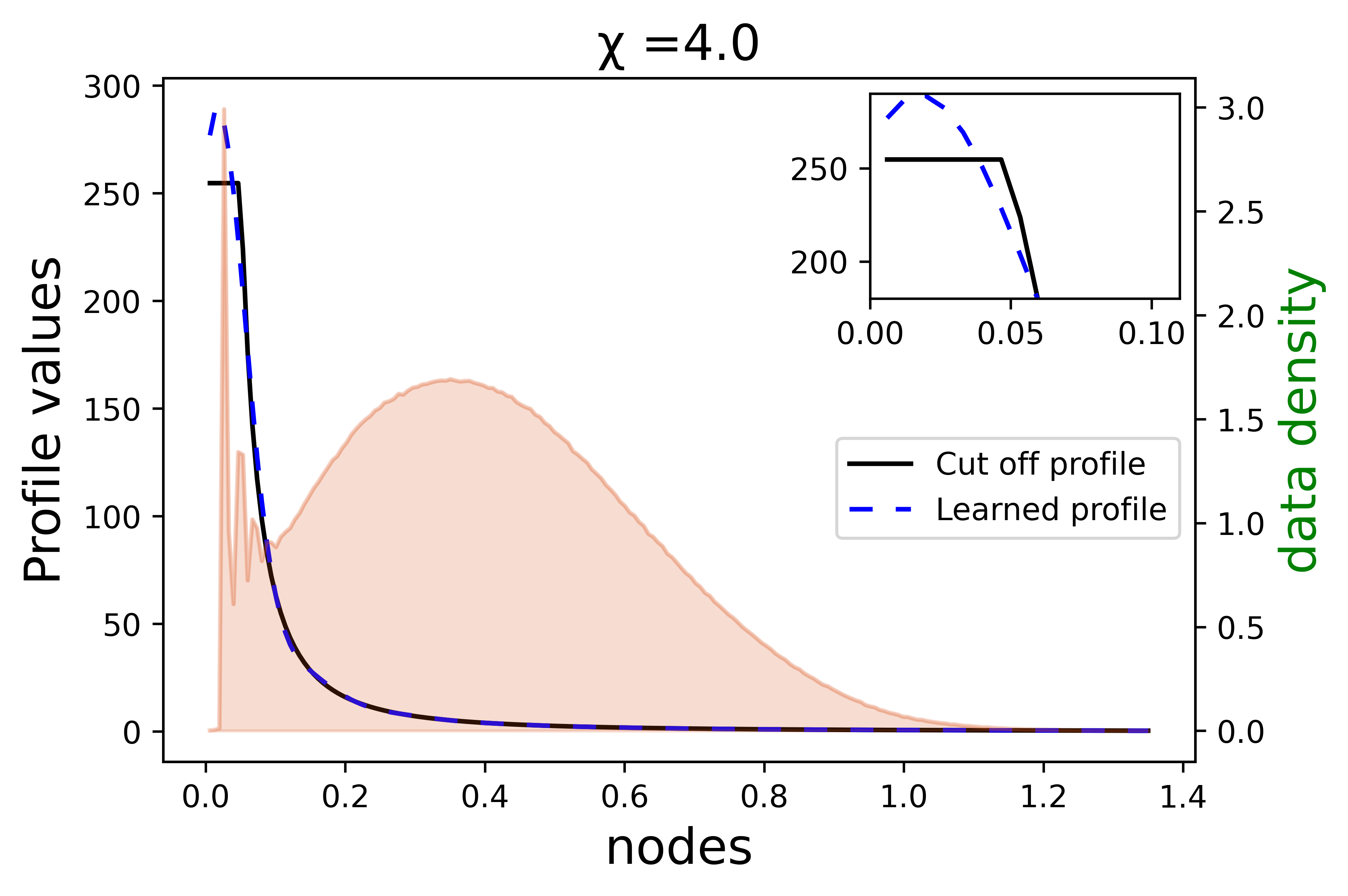}
        \caption{}
    \end{subfigure}
    \begin{subfigure}{0.6\linewidth}
        \centering
        \includegraphics[width=\linewidth]{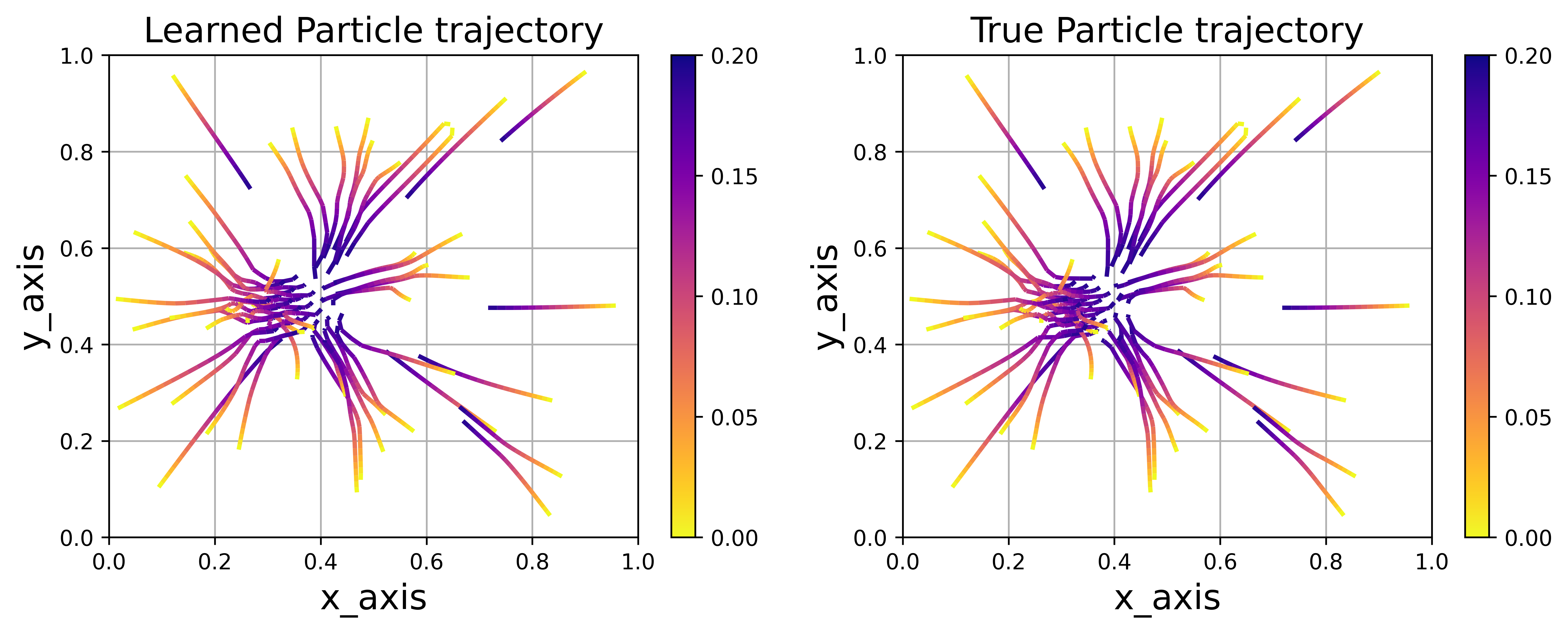}
        \caption{}
    \end{subfigure}
    \caption{Comparison of the learned and regularized profile functions in the two-dimensional Keller–Segel model, obtained using 20 uniform knot points and a truncation parameter of $r_c = 0.05$. The left, middle, and bottom subfigures correspond to initial weight parameters $\chi = 1.0$, $2.0$, and $4.0$, respectively. For each $\chi$ value, the corresponding particle trajectories—learned (left) and true (right)—are also shown.}  
    \label{fig:2D_kernel_uniform_knot}
\end{figure}

\begin{figure}[h]
    \centering
    \captionsetup{width=\linewidth}
    \begin{subfigure}{0.45\linewidth}
        \centering
        \includegraphics[width=\linewidth]{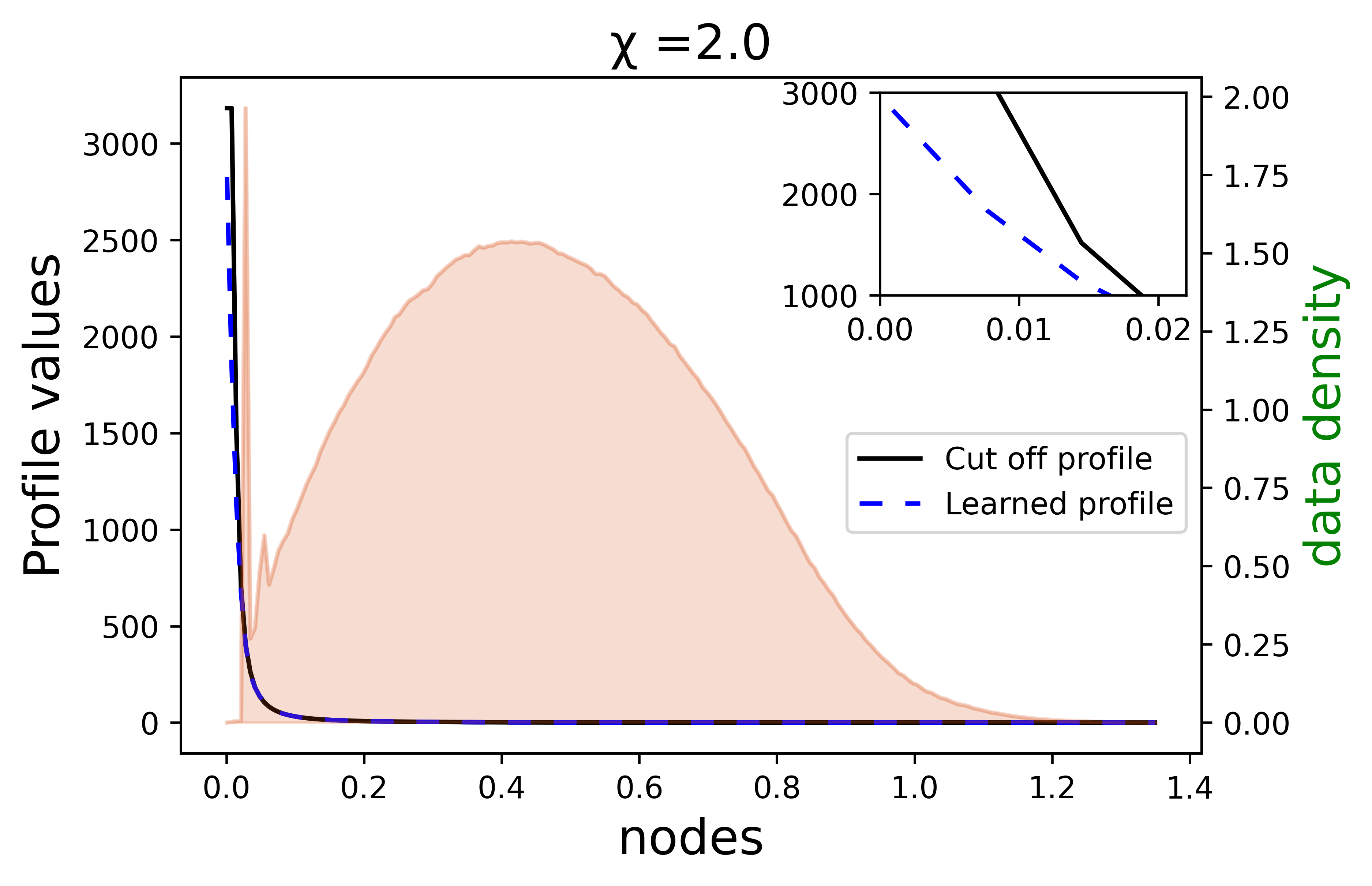}
        \caption{}
    \end{subfigure}
    \hfill
    \begin{subfigure}{0.45\linewidth}
        \centering
        \includegraphics[width=\linewidth]{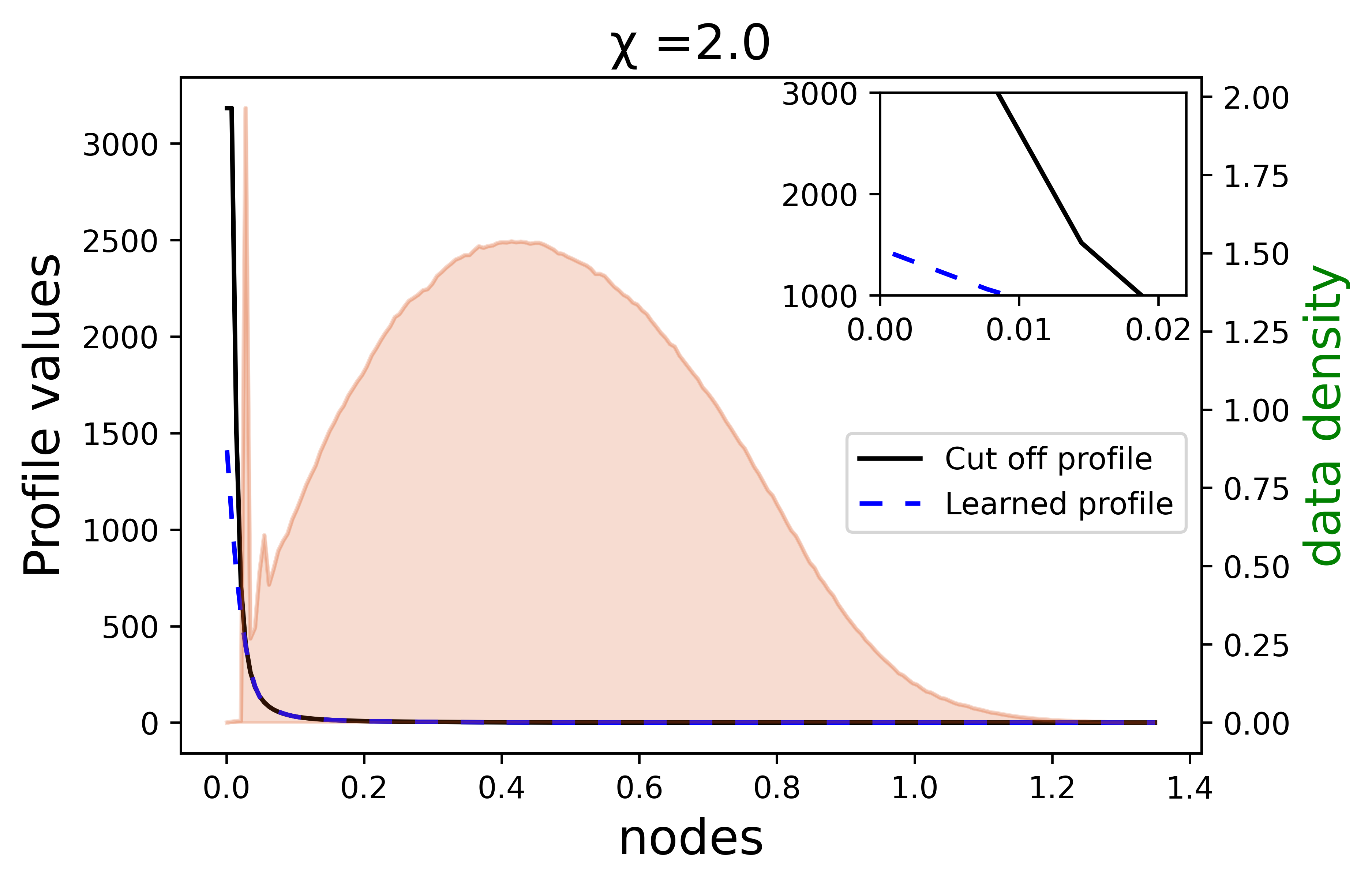}
        \caption{}
    \end{subfigure}    
    \caption{Comparison of the learned profile functions with the regularized profile functions in the two-dimensional Keller–Segel model using 22 adaptive (left) and uniform (right) knot points, with $r_c = 0.01$ and $\omega = 2.0$. Insets provide a detailed view of the profile near the origin.}
    \label{fig:2D_kernel_adaptive_knot}
\end{figure}

In the two-dimensional case, we examine the classical Keller–Segel model, where the Newtonian potential, as defined in \eqref{KS_kernel}, serves as the kernel. To evaluate the accuracy of the learned profile function using uniform knot points, we select the parameters $\chi = 1.0$, $2.0$, and $4.0$. For the adaptive knot points, we set $\chi = 2.0$. The truncation parameter $r_c$ is set to $0.05$ for the uniform knot case and $0.01$ for the adaptive knot case. \\ 
\indent Figure~\ref{fig:2D_kernel_uniform_knot} compares the learned profile functions $\phi^K$ and the corresponding reconstructed particle trajectories with their regularized counterparts in the $2$-dim Keller–Segel model, obtained using 20 uniform knot points. The data density is high near the origin and drops sharply as the pairwise distance slightly increases. After this initial drop, the density rises again at moderate distances before gradually decaying as the distance continues to increase. From the profile plots, we observe that the learned profiles match well with the regularized profiles before the cut-off point $r_c = 0.05$, even in regions with lower data density. For $r < 0.05$, the use of cubic B-splines as basis functions—continuously differentiable up to the third order—ensures that the profile remains smooth, in contrast to the flat region near the cut-off.
In the trajectory plots, the reconstructed particle trajectories—obtained using the learned profile functions—closely reproduce the regularized trajectories. The results indicate that the aggregation dynamics are well captured by the approximate profile function $\phi^K$ in the two-dimensional case. Moreover, increasing the parameter values leads to a more pronounced aggregation phenomenon.
 \\
\indent Figure~\ref{fig:2D_kernel_adaptive_knot} compares the learned profile functions obtained using uniform and adaptive knot points, with the number of knots set to 22 for both cases. The results clearly indicate that the profile obtained with adaptive knot points aligns significantly more closely with the regularized function than that obtained with uniform knots, demonstrating the effectiveness of adaptive knot placement in the two-dimensional setting.  

\indent Table~\ref{tab:2D} reports the numerical errors in reconstructed trajectories and learned profiles for the two-dimensional case. Compared with the one-dimensional results, here we use a larger cutoff value, $r_c = 0.05$, which helps reduce the data errors. The trajectory errors remain small across all values of the sensitivity parameter $\chi$, while the profile errors increase slightly as $\chi$ grows, ranging from $4.77\mathrm{e}{-2}$ at $\chi = 1.0$ to $6.96\mathrm{e}{-2}$ at $\chi = 4.0$. Overall, as in the one-dimensional case, the trajectory reconstruction is robust, and the singularity of the kernel primarily affects the accuracy of the profile reconstruction. 
\begin{table}[h!]
  \centering 
  \captionsetup{width=\linewidth}  
  \renewcommand{\arraystretch}{1.3} 
  \setlength{\tabcolsep}{6pt} 
  \begin{minipage}{\columnwidth}
  \centering
    
  \begin{tabular}{ |c | c c c |}
    \hline
    $\chi$ & $1.0$ & $2.0$ & $4.0$ \\ \hline
    $Err_{\text{traj}}^{\text{rel}}$ & $5.89\mathrm{e}{-3}$ & $3.16\mathrm{e}{-4}$ & $5.96\mathrm{e}{-4}$ \\ \hline
    $Err_{\phi}^{\text{rel}}$ & $4.77\mathrm{e}{-2}$ & $5.96\mathrm{e}{-2}$ & $6.96\mathrm{e}{-2}$ \\ \hline
  \end{tabular}
  \caption{Numerical errors in reconstructed trajectories and learned profiles (2D case) obtained from the learned profile function using 20 uniform knots and $r_{c} = 0.05$.}
  \label{tab:2D}
  \end{minipage}
\end{table} 
\subsection{The Three-Dimensional Case}
\begin{figure}[h]
    \centering 
    \captionsetup{width=\linewidth}
    \begin{subfigure}{0.37\linewidth}
        \centering
        \includegraphics[width=\linewidth]{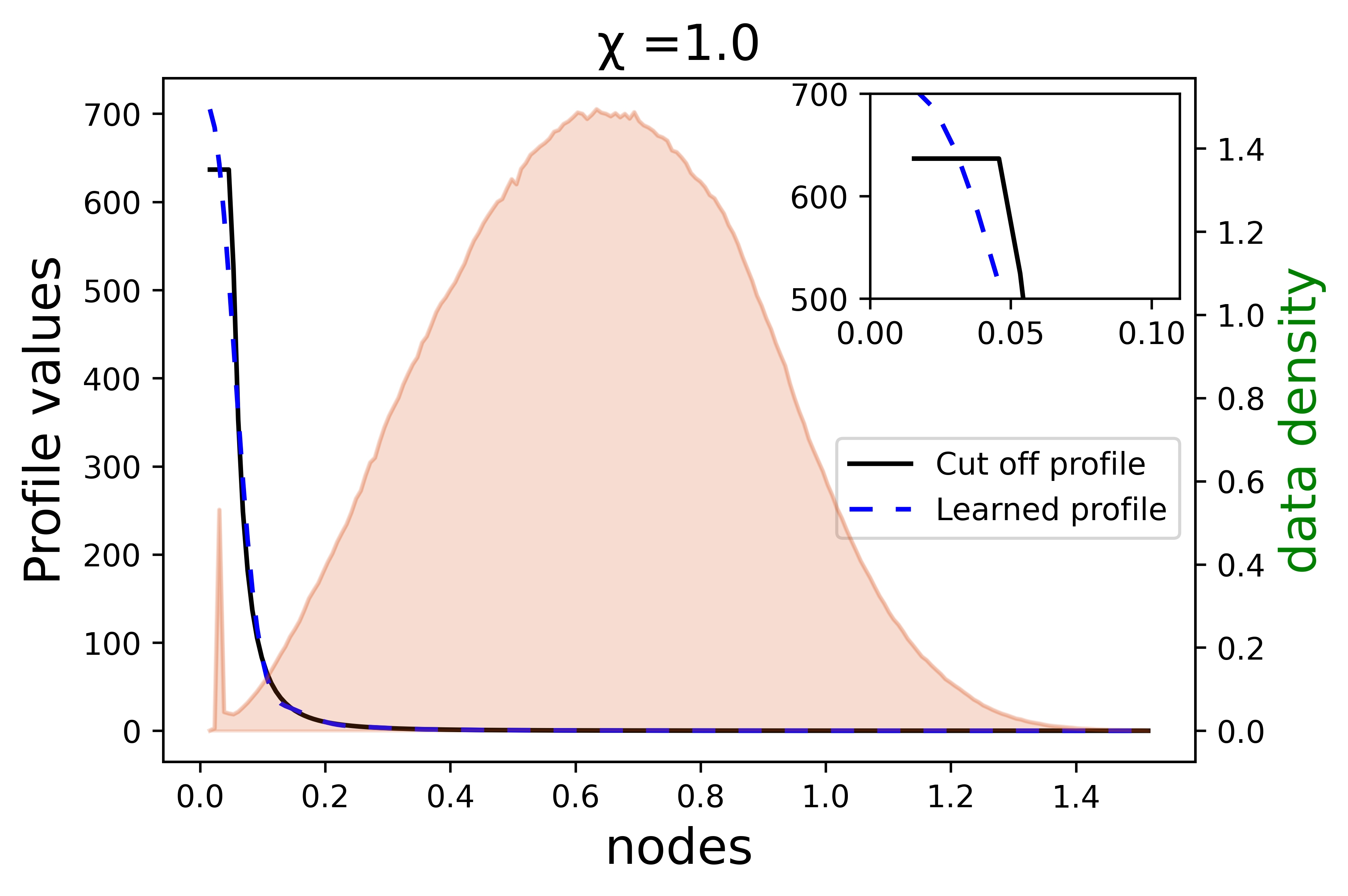}
        \caption{}
    \end{subfigure}
    \begin{subfigure}{0.62\linewidth}
        \centering
        \includegraphics[width=\linewidth]{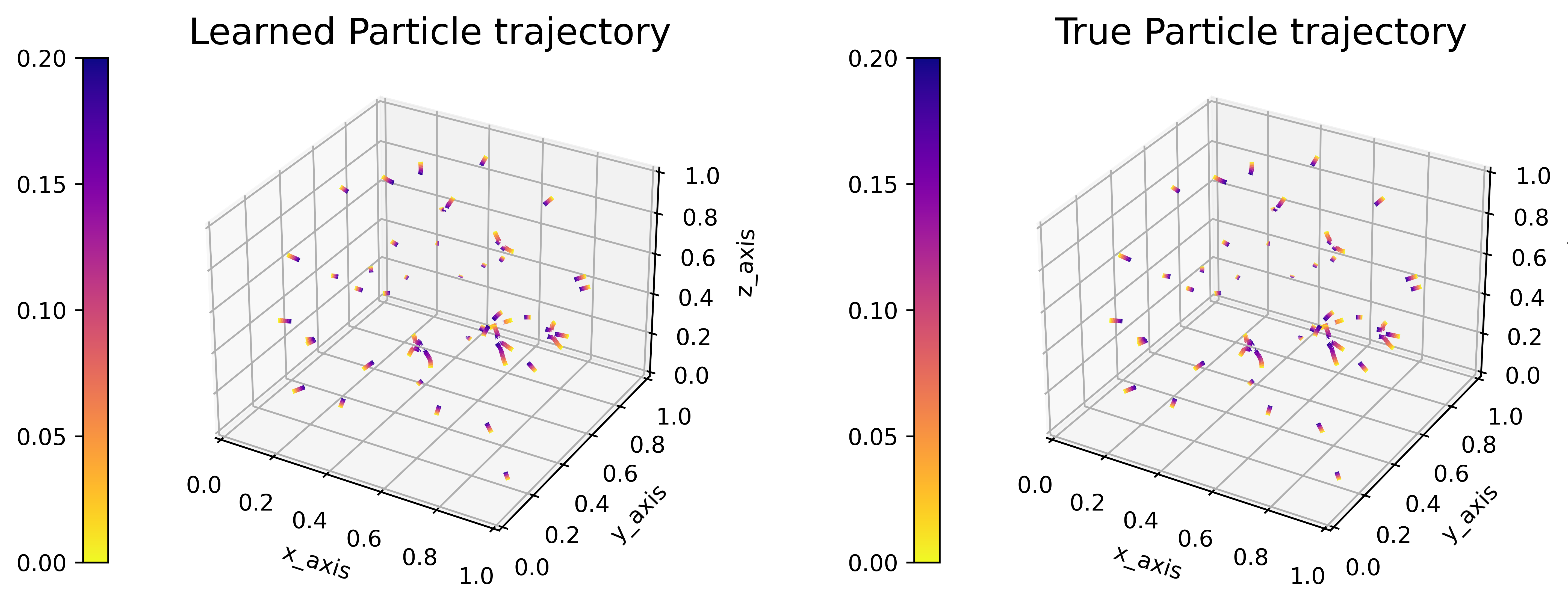}
        \caption{}
    \end{subfigure} 
    \begin{subfigure}{0.37\linewidth}
        \centering
        \includegraphics[width=\linewidth]{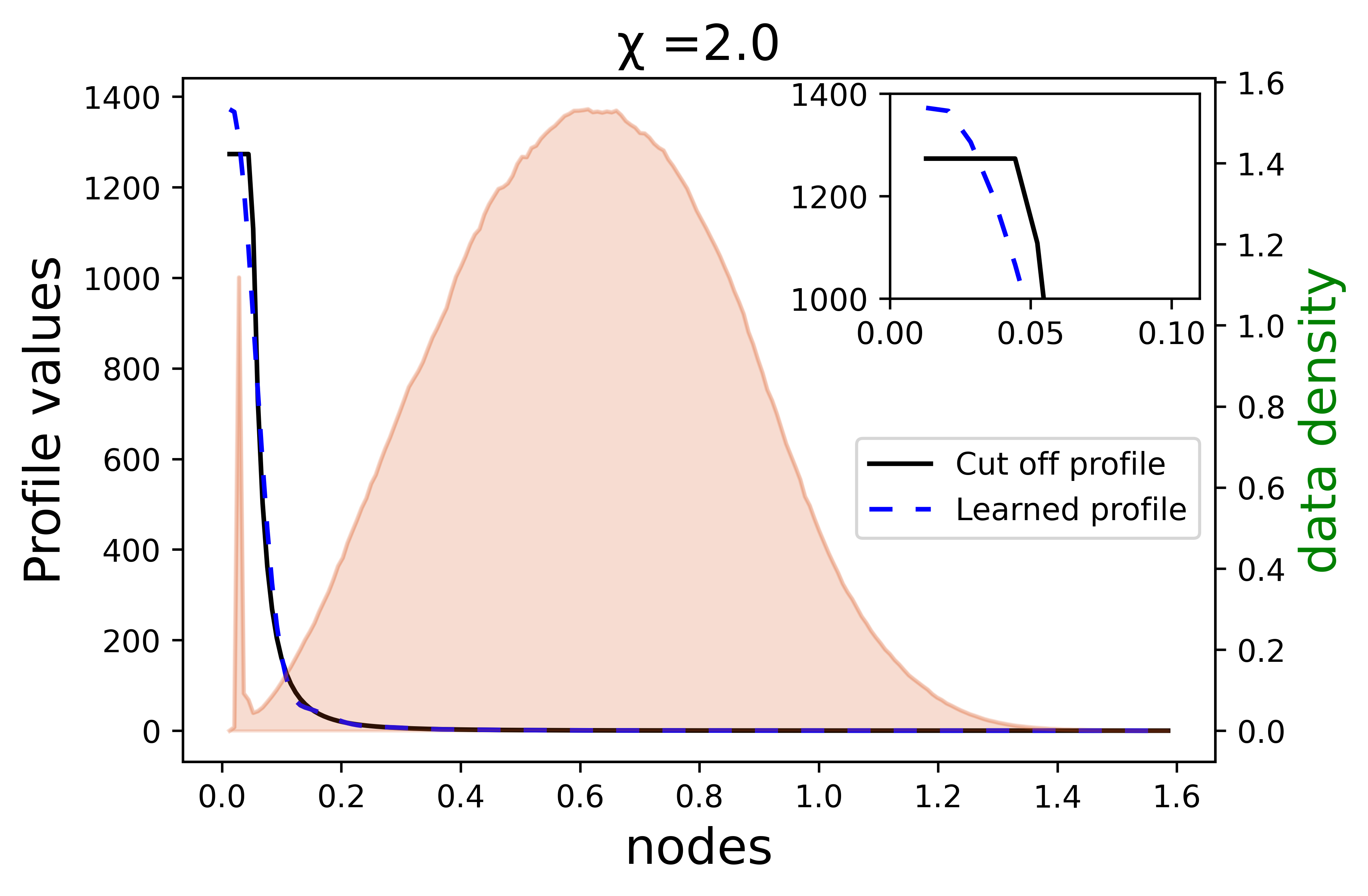}
        \caption{}
    \end{subfigure}
    \begin{subfigure}{0.62\linewidth}
        \centering
        \includegraphics[width=\linewidth]{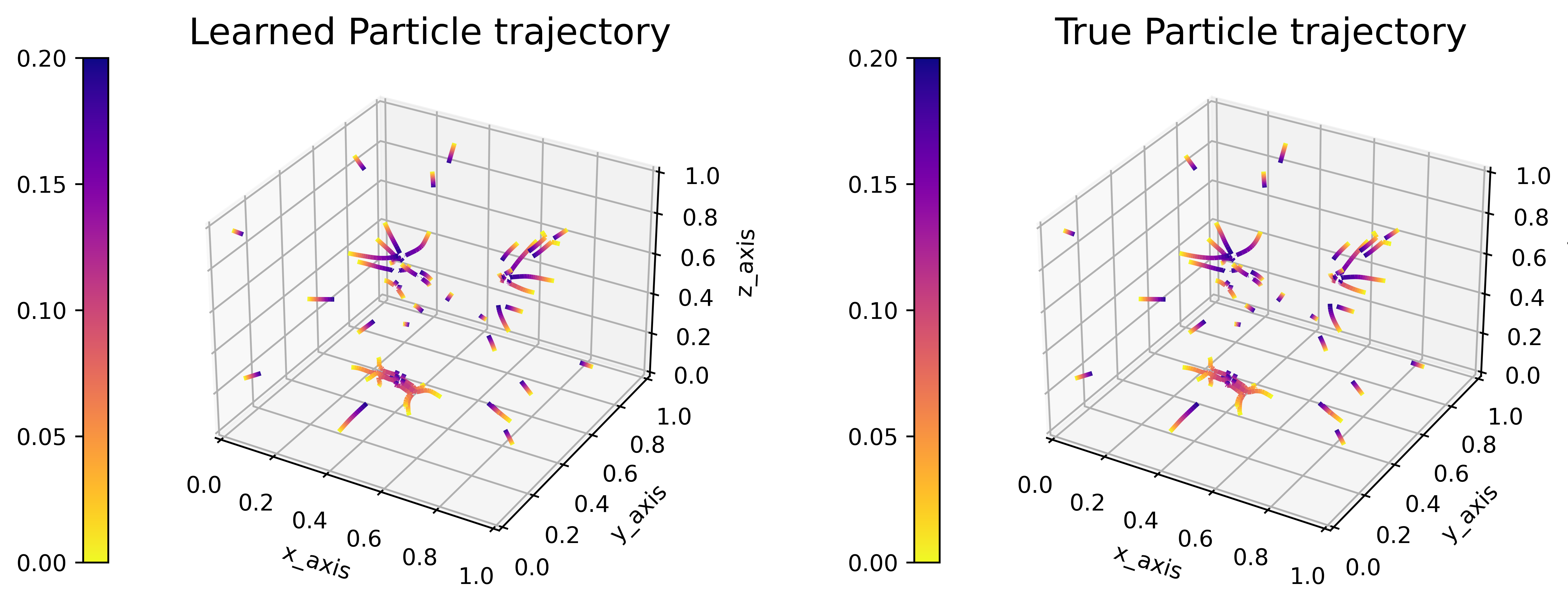}
        \caption{}
    \end{subfigure}
    \begin{subfigure}{0.37\linewidth}
        \centering
        \includegraphics[width=\linewidth]{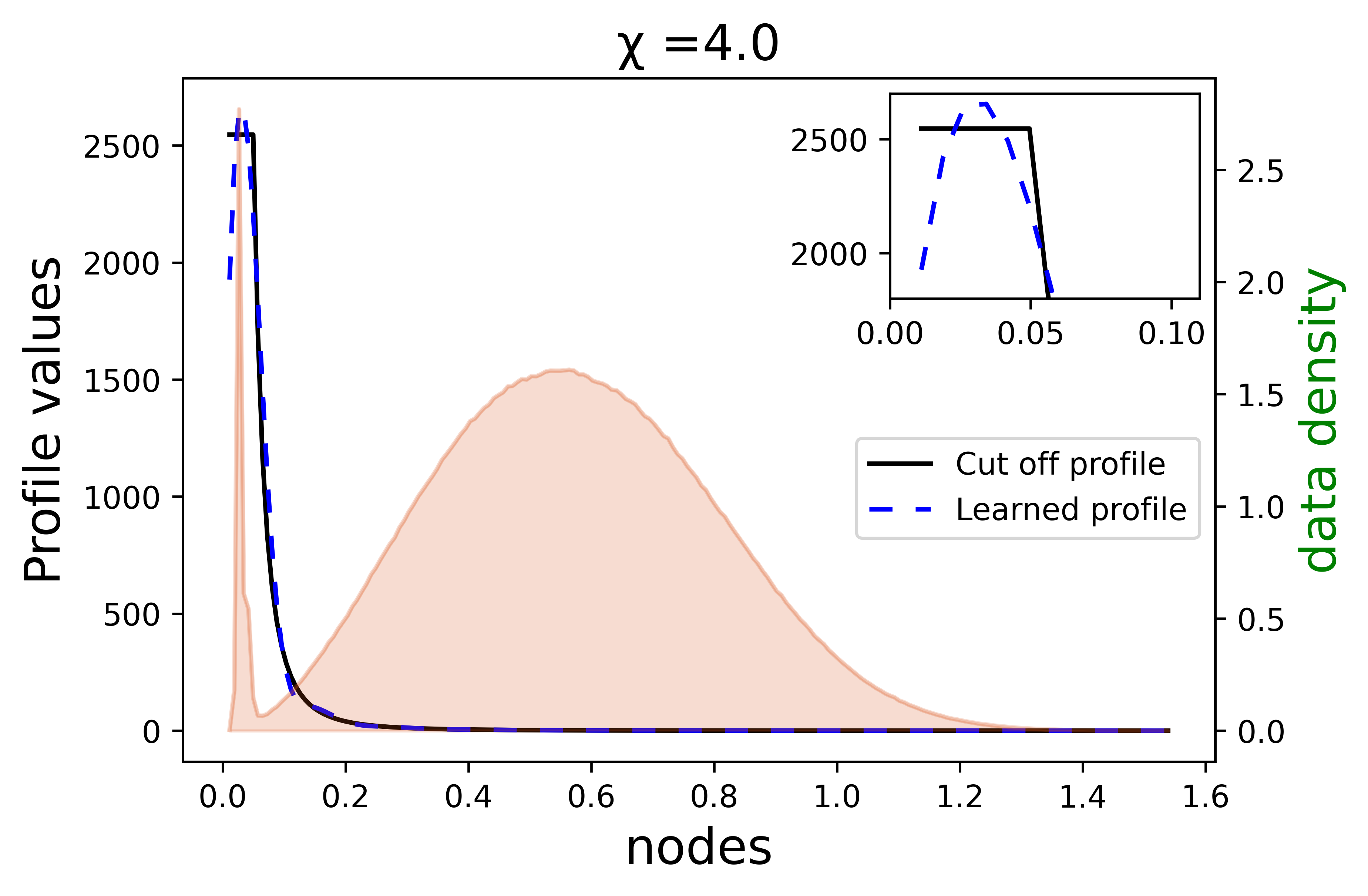}
        \caption{}
    \end{subfigure}
    \begin{subfigure}{0.62\linewidth}
        \centering
        \includegraphics[width=\linewidth]{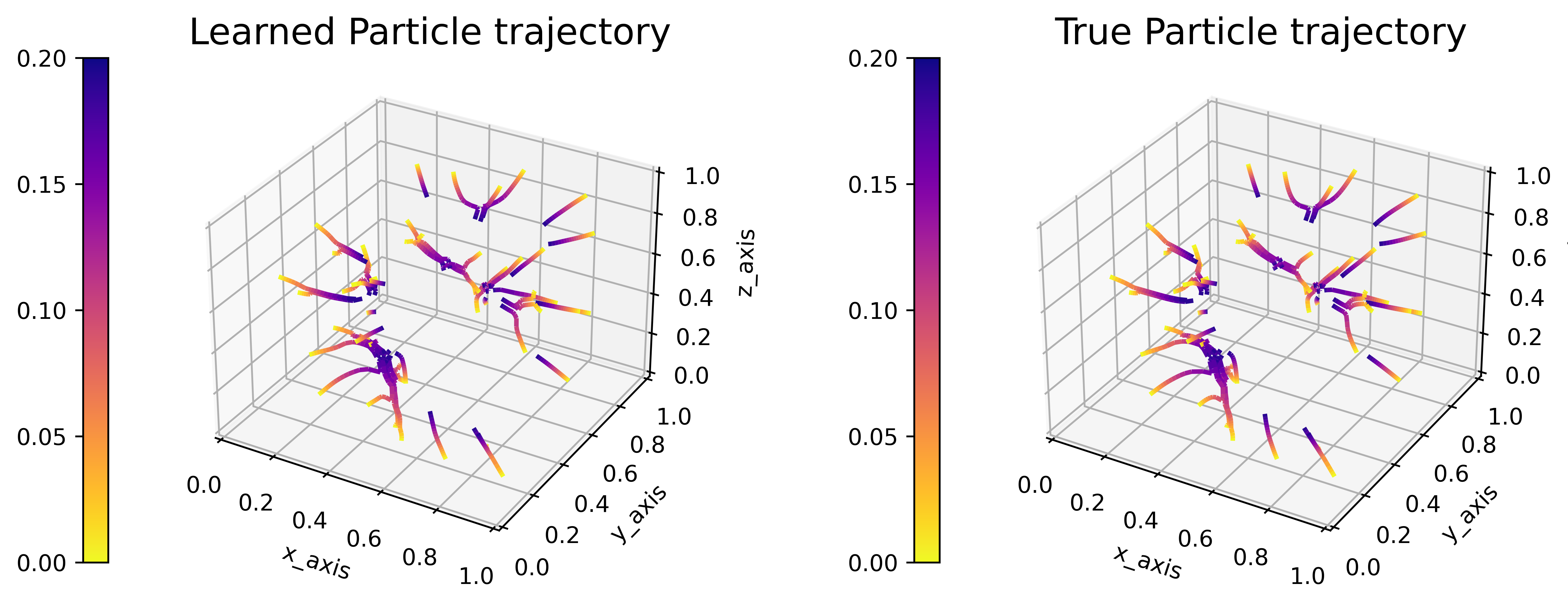}
        \caption{}
    \end{subfigure}
    \caption{Comparison of the learned profile functions with the regularized profile functions in the three-dimensional modified Keller–Segel model, using 25 uniform knot points and a truncation parameter of $r_c = 0.05$. The left, middle, and bottom subfigures correspond to $\chi = 1.0$, $2.0$, and $4.0$, respectively. The corresponding three-dimensional particle trajectories—learned (left) and true (right)—are also shown, using the same parameter settings. The colorbar represents time.}
    \label{fig:3D_kernel_uniform_knot}
\end{figure}

In the three-dimensional case, we assess the accuracy of the learned profile function using uniform knot points with the same initial weight parameters as in the two-dimensional setting: $\chi = 1.0$, $2.0$, and $4.0$. For the adaptive knot points, we again use $\chi = 2.0$. The truncation parameter is set to $r_c = 0.05$ for the uniform knot case and $r_c = 0.01$ for the adaptive knot case.
Figure~\ref{fig:3D_kernel_uniform_knot} compares the learned profile functions $\phi^K$ and the corresponding reconstructed particle trajectories with their regularized counterparts in the three-dimensional modified Keller–Segel model, obtained using 25 uniform knot points. The data density exhibits a similar trend to the two-dimensional case but decreases even more sharply as the pairwise distance slightly increases and remains low over a longer range. From the profile plots, we observe that the learned profiles closely match the regularized profiles up to the cut-off point $r_c = 0.05$, even in regions with sparse data. For $r < 0.05$, the use of cubic B-splines as basis functions—continuously differentiable up to the third order—ensures that the profile remains smooth, in contrast to the flat segment near the cut-off.
In the trajectory plots, the reconstructed particle trajectories—obtained using the learned profile functions—closely reproduce the regularized trajectories. Similar to the two-dimensional case, the results indicate that the aggregation dynamics are well captured by our approximate profile function $\phi^K$ in three dimensions. Moreover, consistent with the 2D observations, increasing the initial weight leads to a more pronounced aggregation phenomenon.

\begin{figure}[ht]
    \centering 
    \captionsetup{width=\linewidth}
    \begin{subfigure}{0.45\linewidth}
        \centering
        \includegraphics[width=\linewidth]{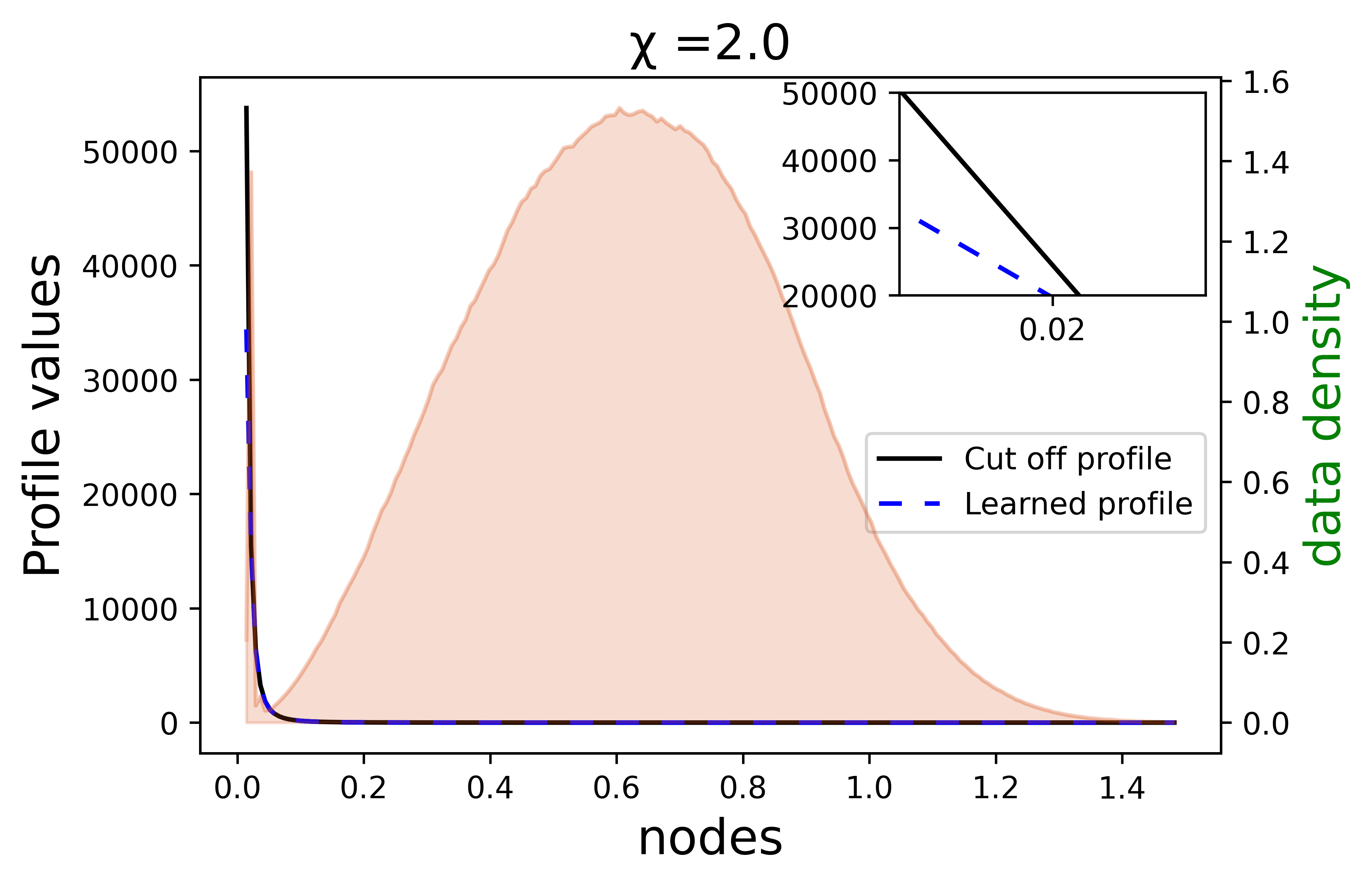}
        \caption{}
    \end{subfigure}
    \hfill
    \begin{subfigure}{0.45\linewidth}
        \centering
        \includegraphics[width=\linewidth]{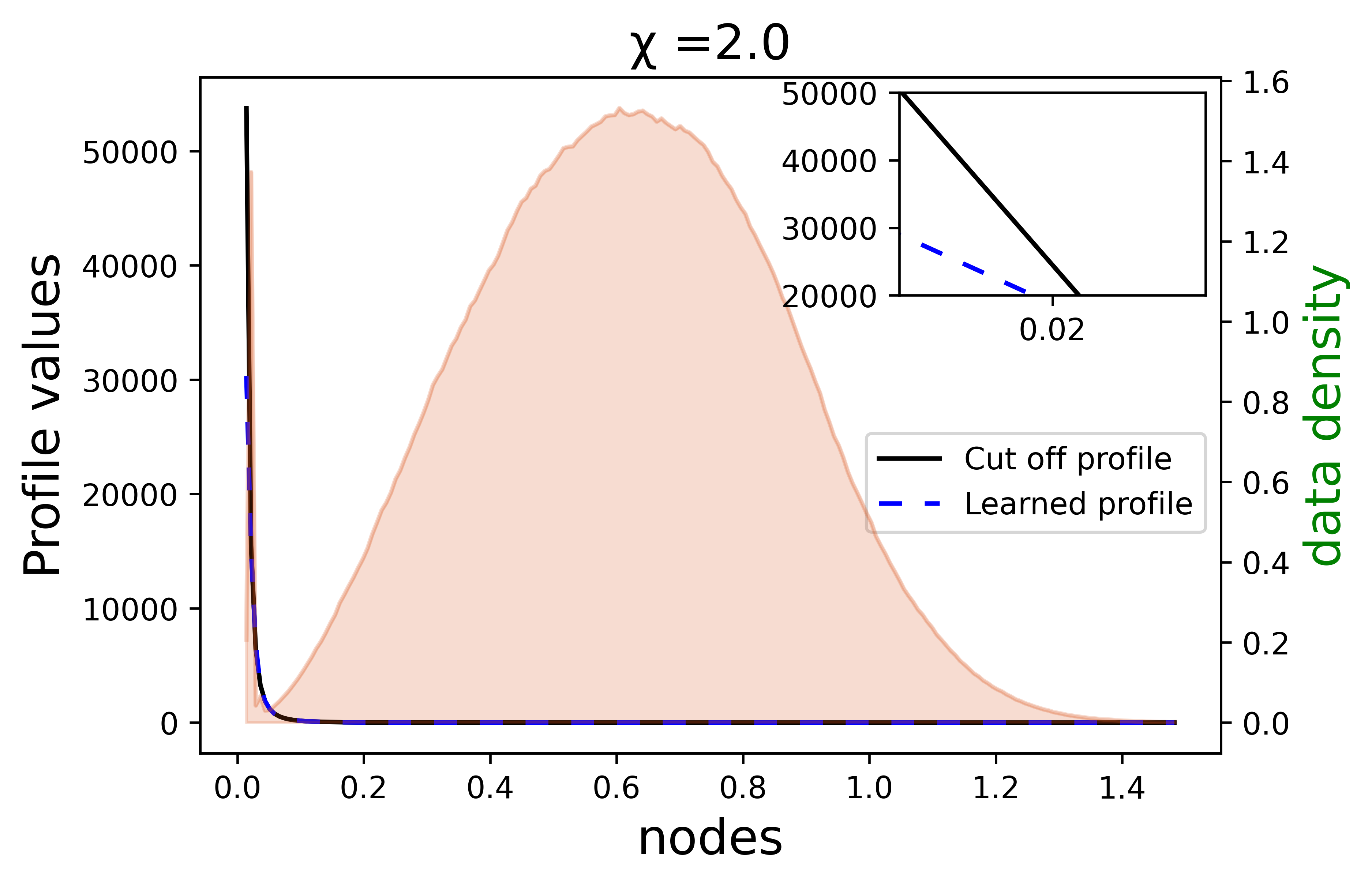}
        \caption{}
    \end{subfigure}
    \caption{Comparison of the learned profile functions with the regularized profile functions in the three-dimensional Keller–Segel model, using 52 adaptive (left) and uniform (right) knot points, with $r_c = 0.01$ and $\omega = 2.0$. Localized behavior of the profile near the origin is illustrated in the insets}
    \label{fig:3D_kernel_adaptive_knot}
\end{figure}

\indent Figure~\ref{fig:3D_kernel_adaptive_knot} compares the learned profile functions obtained using uniform and adaptive knot points in the three-dimensional case, with the number of knots set to 52 for both cases. The adaptive knot point case continues to outperform the uniform case. It is also worth noting that the kernel values are significantly larger in scale compared to the two-dimensional case.
 \\ 

\indent Table~\ref{tab:3D} reports the numerical errors in reconstructed trajectories and learned profiles for the three-dimensional case. The learned profile function is obtained using 25 uniform knots with a cutoff value $r_c = 0.05$, as in the two-dimensional case. We observe that the trajectory errors remain small, ranging from $6.41\mathrm{e}{-4}$ at $\chi = 1.0$ to $3.21\mathrm{e}{-3}$ at $\chi = 4.0$, while the profile errors stay moderate and, unlike in the lower-dimensional cases, do not increase with $\chi$, ranging from $0.105$ to $0.072$. Compared with the two-dimensional case, the larger numerical error can be attributed to the higher values of the profile function near the singularity.

\begin{table}[h!]
  \centering 
  \captionsetup{width=\linewidth}  
  \renewcommand{\arraystretch}{1.3} 
  \setlength{\tabcolsep}{6pt} 
  \begin{minipage}{\columnwidth}
  \centering

  \begin{tabular}{ |c | c c c |}
    \hline
    $\chi$ & $1.0$ & $2.0$ & $4.0$ \\ \hline
    $Err_{\text{traj}}^{\text{rel}}$ & $6.41\mathrm{e}{-4}$ & $1.43 \mathrm{e}{-3}$ & $3.21\mathrm{e}{-3}$ \\ \hline
    $Err_{\phi}^{\text{rel}}$ & $0.105$ & $0.103$ & $0.072$ \\ \hline
  \end{tabular}
  \caption{Numerical errors in reconstructed trajectories and learned profiles (3D case) obtained from the learned profile function using 25 uniform knots and $r_{c} = 0.05$.}
  \label{tab:3D}
  \end{minipage}
\end{table}

\subsection{The Four-Dimensional Case}  
\begin{figure}[ht]
    \centering
    \captionsetup{width=\linewidth}
    \begin{subfigure}{0.5 \linewidth}
        \centering
        \includegraphics[width=\linewidth]{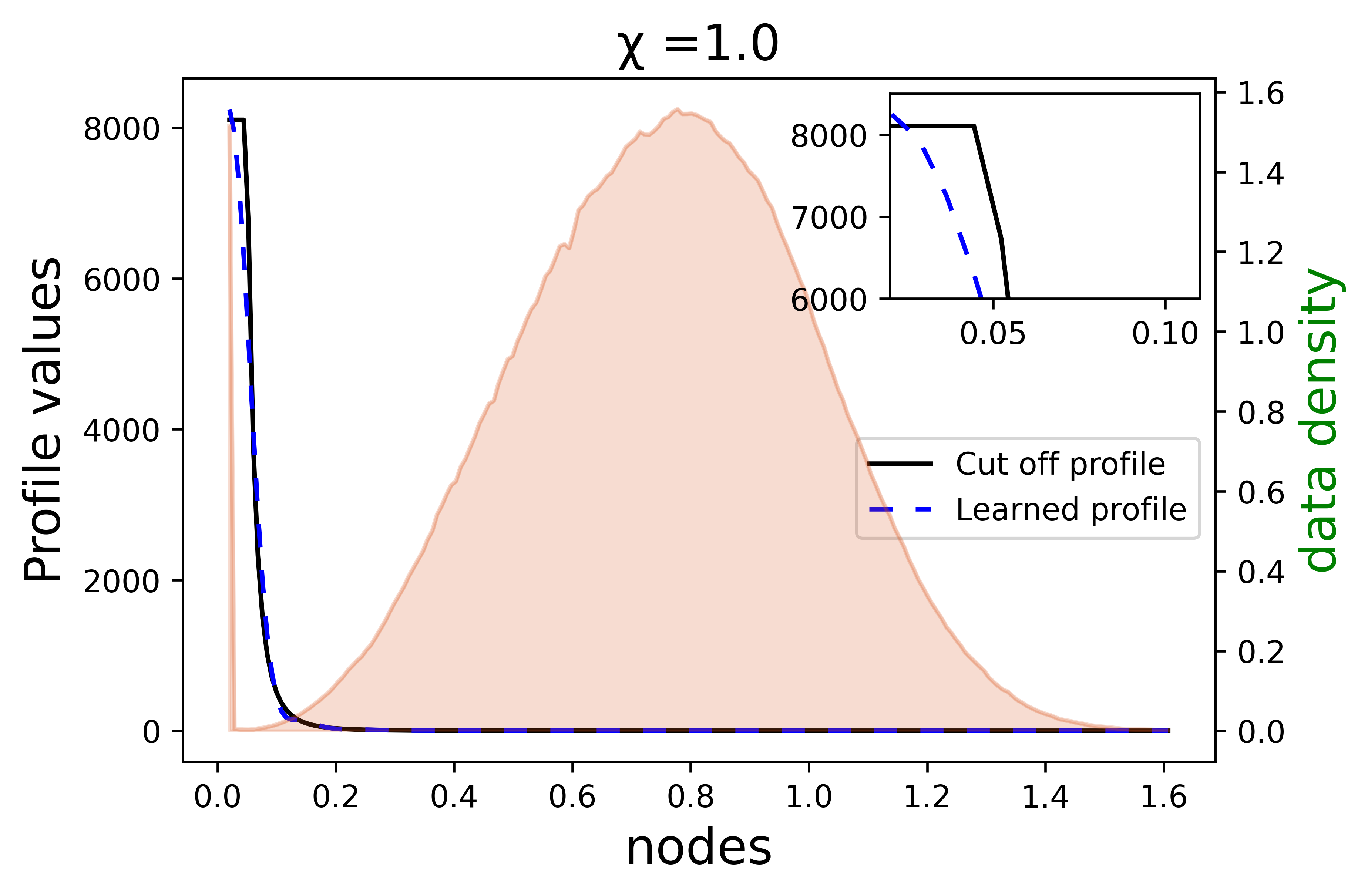}
    \end{subfigure}
    \caption{Comparison of the learned profile functions with the regularized profile functions in the four-dimensional Keller–Segel model, using 30 uniform knot points, a truncation parameter of $r_c = 0.05$, and an initial weight of $\omega = 1.0$. The inset provides a zoomed-in view of the profile near the origin.}
    \label{fig:4D_kernel_uniform_knot}
\end{figure}

In the four-dimensional case, we restrict our attention to the initial weight $\omega = 1.0$ due to the model’s tendency to exhibit rapid particle aggregation. Larger weights result in earlier blow-up, making it infeasible to obtain reliable data for training. For the uniform knot case, the truncation parameter $r_c$ is set to $0.05$. 
Figure~\ref{fig:4D_kernel_uniform_knot} shows a comparison between the learned profile function $\phi^K$ and the regularized profile function $\tilde{\phi}$, using 30 uniform knot points. Similar to the three-dimensional case, the data density decays rapidly with increasing pairwise distance and remains low over a longer range. The learned profile closely follows the regularized profile up to the cut-off point $r_c = 0.05$, with good agreement even in regions with sparse data, and maintains accuracy near the origin. 

\indent For the four-dimensional case with $\chi = 1.0$, the numerical errors are $Err_{\text{traj}}^{\text{rel}} = 1.44 \mathrm{e}{-3}$ for the particle trajectories and $Err_{\phi}^{\text{rel}} = 6.88 \mathrm{e}{-2}$ for the learned profile function using 30 uniform knots and $r_c = 0.05$.

\pagebreak 
\subsection{The Two-Dimensional Stochastic case} 
In the stochastic case, we consider the Keller–Segel model perturbed by Brownian noise, where the Newtonian potential, as defined in \eqref{KS_kernel}, again serves as the kernel. The purpose of this test is to evaluate the robustness of the learned profile function in the presence of stochastic particle trajectories. As in the deterministic case, we select the sensitivity parameters $\chi =$ 1.0, 2.0, and 4.0, while for the adaptive knot points we set the configuration corresponding to $\chi = 2.0$. We choose the regularization parameter $\epsilon = 0.01$.
\begin{figure}[ht]
\centering 
\captionsetup{width=\linewidth}
 \begin{subfigure}{0.37\linewidth}
        \centering
        \includegraphics[width=\linewidth]{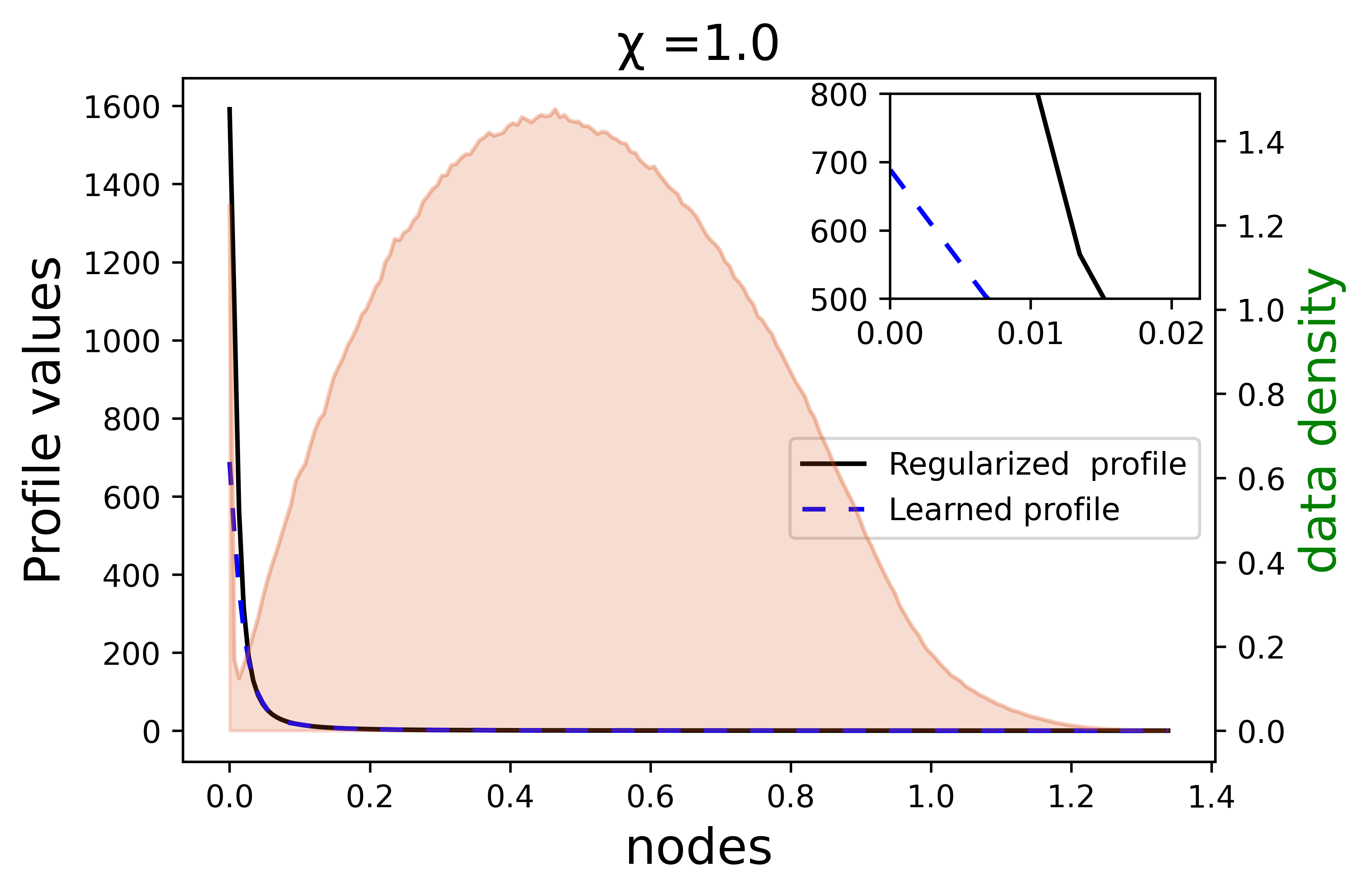}
        \caption{}
    \end{subfigure}
    \begin{subfigure}{0.6\linewidth}
        \centering
        \includegraphics[width=\linewidth]{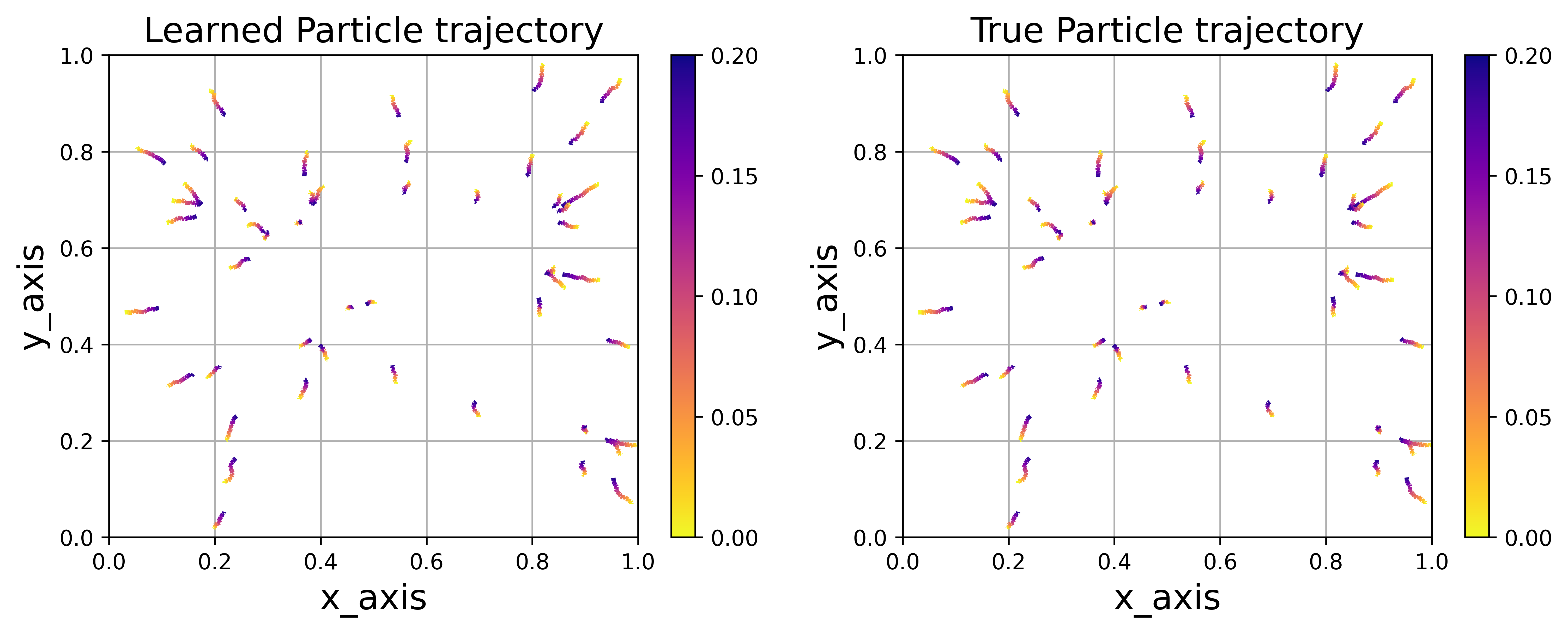}
        \caption{}
    \end{subfigure}
    \begin{subfigure}{0.37\linewidth}
        \centering
        \includegraphics[width=\linewidth]{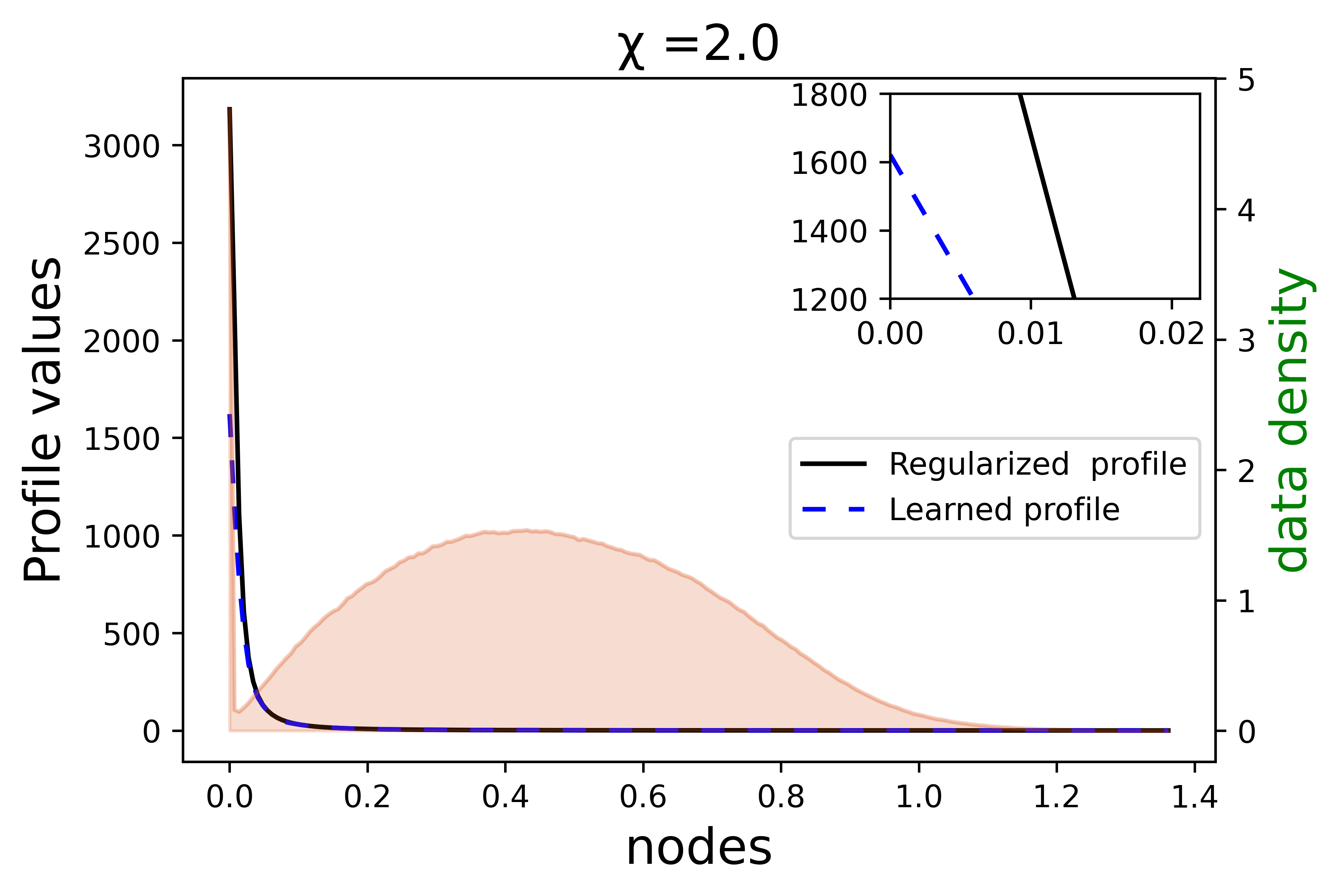}
        \caption{}
    \end{subfigure}
    \begin{subfigure}{0.6\linewidth}
        \centering
        \includegraphics[width=\linewidth]{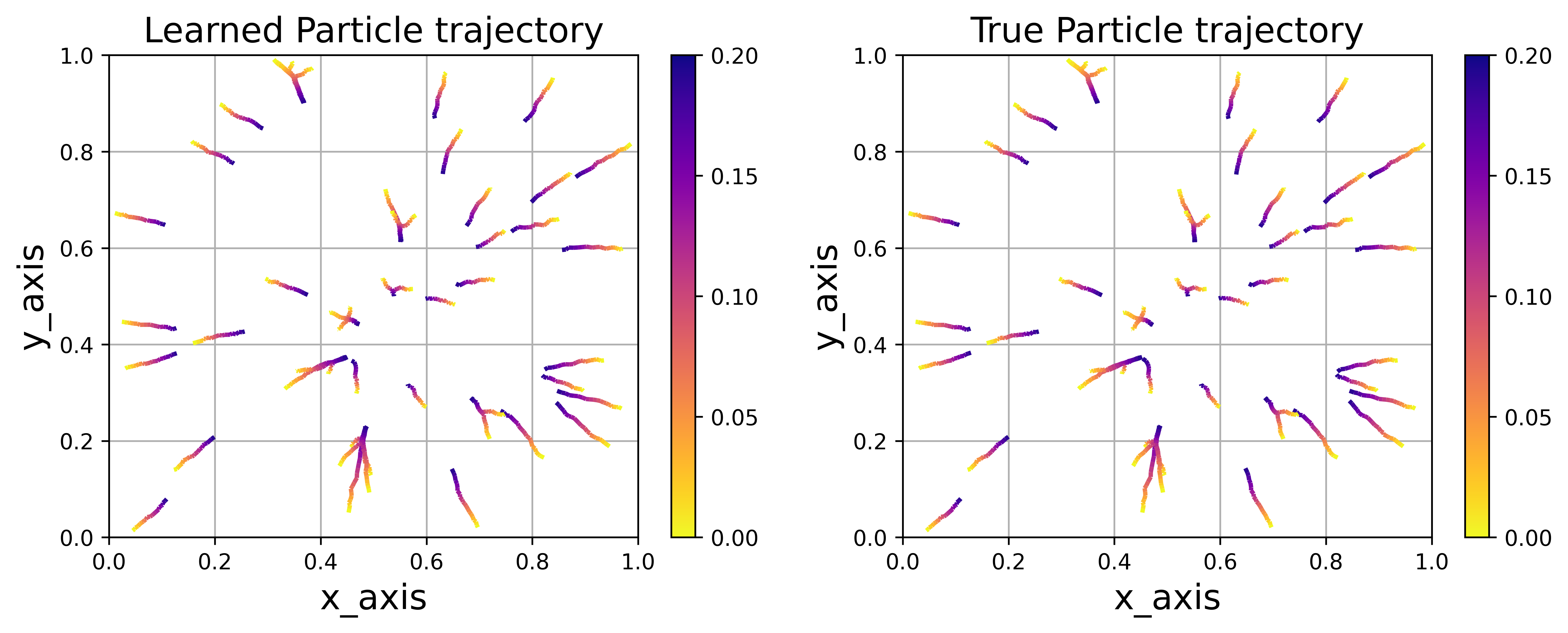}
        \caption{}
    \end{subfigure}
    \begin{subfigure}{0.37\linewidth}
        \centering
        \includegraphics[width=\linewidth]{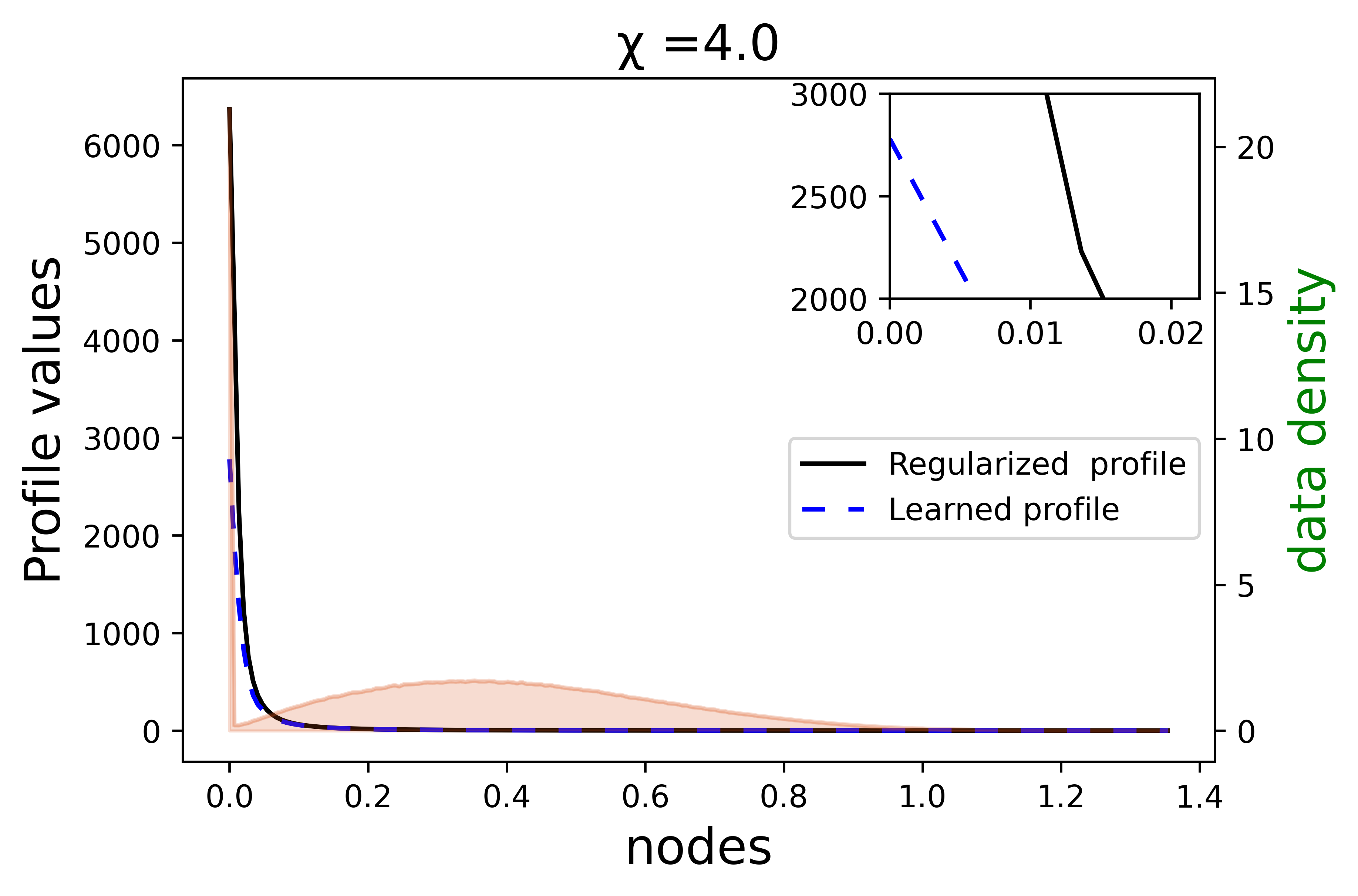}
        \caption{}
    \end{subfigure}
    \begin{subfigure}{0.6\linewidth}
        \centering
        \includegraphics[width=\linewidth]{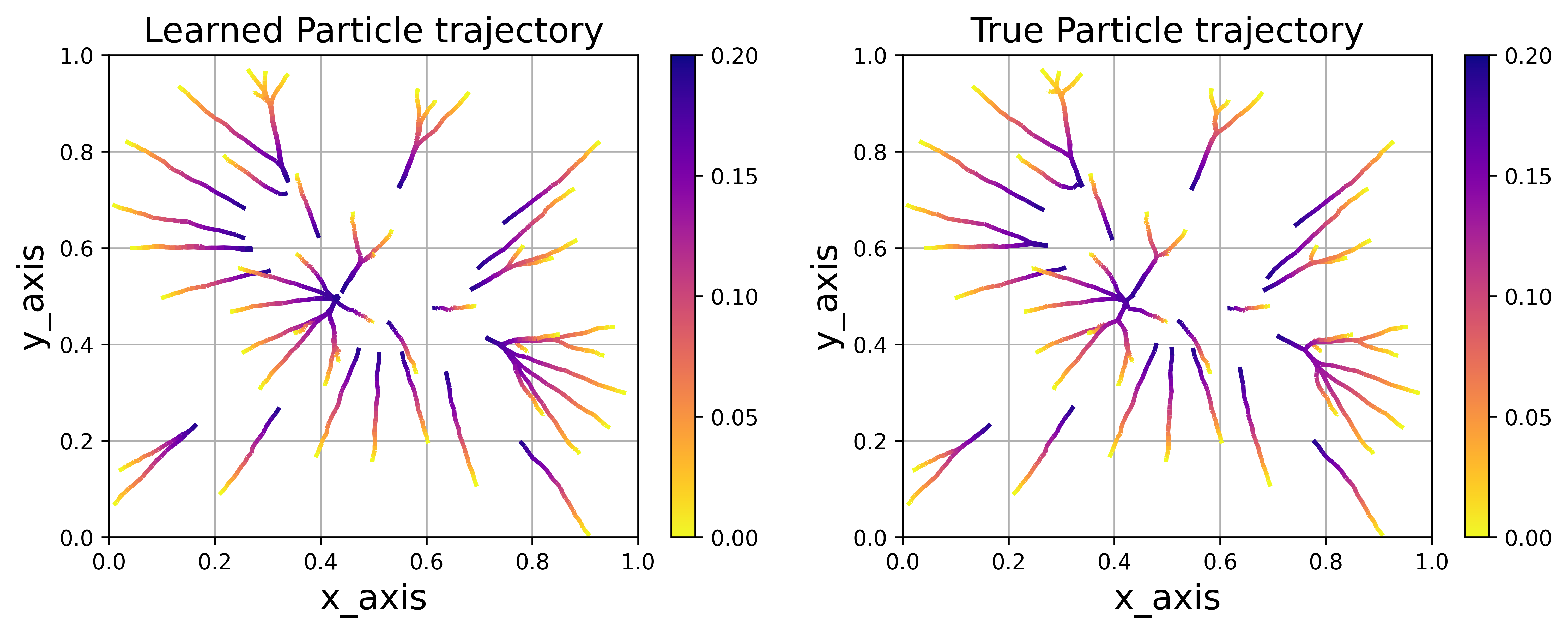}
        \caption{}
    \end{subfigure}
\caption{Comparison of the learned and regularized profile functions in the two-dimensional stochastic Keller–Segel model, obtained using 30 uniform knot points and a regularization parameter of $\epsilon = 0.01$. The left, middle, and bottom subfigures correspond to initial weight parameters $\chi = 1.0$, $2.0$, and $4.0$, respectively. For each $\chi$ value, the corresponding particle trajectories—learned (left) and true (right)—are also shown.}
\label{fig:2D_SDE_kernel}
\end{figure}

Figure~\ref{fig:2D_SDE_kernel} presents a comparison between the learned profile function $\phi^K$ and the regularized profile function $\phi^{\epsilon}$, computed using 30 uniform knot points, along with the corresponding particle trajectories. In the stochastic case, the data density attains higher values near the origin as $\chi$ increases, and the learned profiles lie below the regularized profiles close to the origin, reflecting the effect of noise. The reconstructed particle trajectories closely reproduce the regularized stochastic trajectories, showing that the aggregation dynamics are accurately captured, with additional small-scale fluctuations introduced by the diffusion term through the parameter $\eta$.

\indent Table~\ref{tab:S2D} reports the numerical errors in reconstructed stochastic particle trajectories and in the learned profile functions for the two-dimensional case. We observe that the trajectory errors remain small across all values of the sensitivity parameter $\chi$, increasing slightly from $4.29\mathrm{e}{-4}$ at $\chi = 1.0$ to $3.35\mathrm{e}{-3}$ at $\chi = 4.0$. In contrast, the profile errors are larger and grow with increasing $\chi$, ranging from $0.233$ to $0.380$, indicating that learning the profile becomes more challenging as the strength of the interaction increases. Overall, while stochastic fluctuations slightly increase the trajectory errors compared with the deterministic case, the reconstruction method remains robust, and the trend in profile error with respect to $\chi$ is consistent with the one-dimensional and deterministic two-dimensional cases.

\begin{table}[h!]
  \centering 
  \captionsetup{width=\linewidth}  
  \renewcommand{\arraystretch}{1.3} 
  \setlength{\tabcolsep}{6pt} 
  \begin{minipage}{\columnwidth}
  \centering

  \begin{tabular}{ |c | c c c |}
    \hline
    $\chi$ & $1.0$ & $2.0$ & $4.0$ \\ \hline
    $Err_{\text{traj}}^{\text{rel}}$ & $4.29\mathrm{e}{-4}$ & $6.92 \mathrm{e}{-4}$ & $3.35\mathrm{e}{-3}$ \\ \hline
    $Err_{\phi}^{\text{rel}}$ & $0.233$ & $0.234$ & $0.380$ \\ \hline
  \end{tabular}
  \caption{Numerical errors of reconstructed stochastic particle trajectories and learned profile functions in two dimensions.}
  \label{tab:S2D}
  \end{minipage}
\end{table} 

\section{Conclusion}\label{sec:conclude}
We have demonstrated the effectiveness of our learning method for identifying the profile function of Keller-Segel models in various physical dimensions.  Our learning method can handle the observation data from deterministic or stochastic dynamics without the prior information of the stochastic noise. Moreover, the profile function has a sharp blow-up near zero; our method is able to capture that behavior with the adaptive learning, providing better accuracy with smaller numbers of grid points. Future directions on multi-physics and multi-species Keller-Segel models are ongoing.

\section*{acknowledgments}
MZ is supported by NSF-CCF-$2225507$ and FY$24$ Ralph E. Powe Junior Faculty Enhancement Award by ORAU. C.L. was partially supported by NSF grants DMS-2410742 and DMS-2118181.
%
\bibliographystyle{siamplain}
\bibliography{bibliography}

\end{document}

%% file: ex_shared.tex
\usepackage{lipsum}
\usepackage{amsfonts}
\usepackage{epstopdf}
\usepackage{booktabs}
\usepackage{amsmath}

\usepackage{caption}
\usepackage[font=footnotesize]{caption}
\usepackage{subcaption}
\usepackage{}

\ifpdf
  \DeclareGraphicsExtensions{.eps,.pdf,.png,.jpg}
\else
  \DeclareGraphicsExtensions{.eps}
\fi

\newsiamremark{hypothesis}{Hypothesis}

\newsiamthm{claim}{Claim}

\headers{Learning Profile Functions in Keller{-}Segel Models}{Chi-An Chen, Chun Liu, Ming Zhong}

\title{Unified Learning of the Profile Function in Discrete Keller-Segel Models\thanks{
\funding{MZ’s research is funded by NSF CCF-AOF-2225507 and the Ralph E. Powe Junior Faculty Enhancement Award from ORAU (FY24). CL’s research is partially supported by NSF grants DMS-2410742 and DMS-2118181.}}}

\author{Chi-An Chen \thanks{Department of Applied Mathematics, Illinois Institute of Technology, Chicago, IL 
  (\email{cchen156@hawk.illinoistech.edu}).}
\and Chun Liu \thanks{Department of Applied Mathematics, Illinois Institute of Technology, Chicago, IL
  (\email{cliu124@illinoistech.edu}).}
\and Ming Zhong \thanks{Department of Mathematics, University of Houston}
    (\email{mzhong3@central.uh.edu}).}

\usepackage{amsopn}

\makeatletter
\newcommand*{\addFileDependency}[1]{
  \typeout{(#1)}
  \@addtofilelist{#1}
  \IfFileExists{#1}{}{\typeout{No file #1.}}
}
\makeatother
